\documentclass[final,3p]{elsarticle}

\usepackage{amssymb}
\usepackage[dvipsnames]{xcolor}
\usepackage{lineno}
\usepackage[fleqn]{amsmath}
\usepackage{amssymb}
\usepackage{etoolbox}
\usepackage{graphicx}
\usepackage{dsfont}
\usepackage{txfonts} 
\usepackage{mathrsfs} 
\usepackage{caption} 
\usepackage{subcaption}
\usepackage{enumitem}       

\usepackage{amsthm} 
\usepackage{cases}
\usepackage{siunitx}
\usepackage{upgreek}

\usepackage{xargs}
\usepackage[colorinlistoftodos,prependcaption]{todonotes}

\newcommandx{\todoModify}[2][1=]{\todo[
   linecolor=red,
   backgroundcolor=red!25,
   bordercolor=red,
   #1]{#2}}
\newcommandx{\todoCheck}[2][1=]{\todo[
    linecolor=blue,
    backgroundcolor=blue!25,
    bordercolor=blue
    ,#1]{#2}}
\newcommandx{\todoAdd}[2][1=]{\todo[
    linecolor=red,
    backgroundcolor=red!25,
    bordercolor=red,
    #1]{#2}}
\newcommandx{\todoDelete}[2][1=]{\todo[
    linecolor=Plum,
    backgroundcolor=Plum!25,
    bordercolor=Plum,#1]{#2}}

\theoremstyle{remark}
\newtheorem*{remark}{Remark}
\newcommand{\toVect}[1]{{\boldsymbol{#1}}} 

\newcommand{\dd}{\text{\,}\mathrm{d}}
\newcommand{\ds}{\mathrm{s}}

\newcommand{\x}{\toVect{x}}
\newcommand{\Displacement}{{\toVect{u}}}

\newcommand{\Traction}{\toVect{t}}
\newcommand{\HighOrderTraction}{\toVect{r}}
\newcommand{\HighOrderJump}{\toVect{j}}

\newcommand{\E}{\toVect{E}}

\newcommand{\Strain}{\toVect{\strain}}
\newcommand{\strain}{\varepsilon}
\newcommand{\CauchyStress}{\toVect{\cauchyStress}}
\newcommand{\cauchyStress}{\widehat{\sigma}}
\newcommand{\HyperStress}{\toVect{\hyperStress}}
\newcommand{\hyperStress}{\widetilde{\sigma}}
\newcommand{\ElectricDisp}{\toVect{\electricDisp}}
\newcommand{\electricDisp}{\widehat{D}}

\newcommand{\Piezo}{\toVect{\piezo}}
\newcommand{\piezo}{e}
\newcommand{\Flexo}{\toVect{\flexo}}
\newcommand{\flexo}{\mu}

\newcommand{\dielec}{\kappa}

\newcommand{\elast}{\varmathbb{C}}

\newcommand{\strGr}{h}

\newcommand{\gradient}{\nabla}

\newcommand{\surfaceGradient}{\nabla^{S}}

\newcommand{\jump}[1]{\left\llbracket#1\right\rrbracket}
\newcommand{\jumpMacro}[1]{\left\langle\!\left\langle#1\right\rangle\!\right\rangle}

\newcommand{\eq}{Eq.~}
\newcommand{\eqs}{Eqs.~}
\newcommand{\fig}{Fig.~}

\newcommand{\ie}{i.e.~}
\newcommand{\eg}{e.g.~}

\newcommand{\um}{~\si{\um}}

\newcommand{\bm}[1]{\text{\boldmath $#1$\unboldmath}}  

\biboptions{}

\begin{document}

\begin{frontmatter}
\title{Computational homogenization of higher-order electro-mechanical materials with built-in generalized periodicity conditions}
\author{J. Barcel\'o-Mercader$^1$, D. Codony$^{1}$, A. Mocci$^1$, I. Arias$^{1,2}$}
\address{$^1$ Laboratori de C\`{a}lcul Num\`{e}ric (LaC\`{a}N), Universitat Polit\`{e}cnica de Catalunya (UPC), Campus Nord UPC-C2, E-08034 Barcelona, Spain}
\address{$^2$ Centre Internacional de M{\`e}todes Num{\`e}rics en Enginyeria (CIMNE), E-08034 Barcelona, Spain}
\begin{abstract} 
We present a formulation for high-order generalized periodicity conditions in the context of a high-order electromechanical theory including flexoelectricity, strain gradient elasticity and gradient dielectricity, with the goal of studying  periodic architected metamaterials. Such theory results in fourth-order governing partial differential equations, and the periodicity  conditions involve continuity across the periodic boundary of primal fields (displacement and electric potential) and their normal derivatives, continuity  of the corresponding dual generalized forces (tractions, double tractions, surface charge density and  double surface charge density). Rather than imposing these conditions numerically as explicit constraints, we develop an approximation space which fulfils generalized periodicity by construction. Our method naturally allows us to impose general macroscopic fields (strains/stresses and electric fields/electric displacements) along arbitrary directions, enabling the characterization  of the material anisotropy.  We apply the proposed method to study periodic architected metamaterials with apparent piezoelectricity. We first verify the method by directly comparing the results with a large periodic structure, then apply it to evaluate the anisotropic apparently piezoelectricity of a geometrically polarized 2D lattice, and finally demonstrate the application of the method in a 3D architected metamaterial. 
\end{abstract}
\begin{keyword}
	Lifshitz-invariant Flexoelectricity  \sep Generalized periodicity  \sep RVE \sep Immersed boundary b-spline
\end{keyword}
\end{frontmatter}

\section{Introduction}

Periodic metamaterials exploit the fact that, by suitably designing the geometry of the period unit cell, the architected material can exhibit physical properties (mechanical, acoustic, electrimagnetic, optical, etc) that the base material lacks \cite{Engetha2006,10.1063/1.3490504,Bertoldi:2017aa,Kadic:2019aa,10.1063/5.0152099}. Computational homogenization allows us to analyze and design such materials accounting for general geometries and material behaviors by considering a representative volume element (RVE) and analyzing its response upon application of macroscopic fields \cite{Geers2010}. These macroscopic fields are introduced through generalized periodic conditions at the boundary of the RVE \cite{hassani1998review}. When the continuum theory describing the materials response is of higher-order, the generalized periodic conditions become more complex to formulate and to implement in a computational homogenization approach \cite{Schmidt2022}. Here, we focus on developing a computational framework accounting for such generalized periodicity conditions in the context of metamaterials for electromechanical transduction made of base materials with higher-order electromechanical couplings. We provide next a more detailed background of this application, which is representative of other problems modeled with higher-order partial differential equations and which can benefit from the approach proposed here. 

The ability of some materials to transduce electric energy into mechanical energy and vice versa is widely exploited in many applications such as sensing \cite{gautschi2002piezoelectric}, actuation \cite{sinha2009piezoelectric}, or energy harvesting \cite{safaei2019review,dagdeviren2016recent}, among others. Most current technologies for electromechanical transduction rely on the piezoelectric effect \cite{guerin2021restriction}, the linear coupling between strain and dielectric polarization. Piezoelectricity is fundamentally restricted by symmetryt and supported only by dielectrics exhibiting non-centrosymmetric ionic or molecular structure. Significant research is devoted to synthesising new piezoelectric materials to overcome the limitations of current ones in terms of brittleness, lead-content, operating temperatures and biocompatibility \cite{jaffe1955properties,haertling1999ferroelectric,jaffe1958piezoelectric,saito2004lead}. In recent years, flexoelectricity engineering has open the path towards generating piezoelectric-like responses in non-piezoelectric dielectrics \cite{EditorialJAP2022}. Flexoelectricity is a coupling between strain-gradient and electric field or conversely, electric field gradient and strain \cite{zubko2013flexoelectric}. Strain gradients break locally the spatial inversion symmetry inducing an electric response in any dielectric material, which can be significant at sub-micron scales. Strain gradients can be achieved by inhomogeneous deformations \cite{Cross2006} or by a suitable material microarchitecture \cite{sharma2007possibility}. In  \cite{mocci2021geometrically}, we proposed a class of geometrically polarized architected dielectrics with apparent piezoelectricity. By considering periodic metamaterials made of non-piezoelectric flexoelectric materials, we showed that (1) geometric polarization (lack of geometric centrosymmetry) of the representative volume element (RVE) and (2) small-scale geometric features subjected to bending are enough to achieve an apparent piezoelectric behavior similar to that of Quartz and lead zirconium titanate (PZT) materials.  The systematic design and computational optimization of the microstructures of such metamaterials requires the efficient simulation periodic unit cells under generalized periodic conditions \cite{mocci2021geometrically,Mawassy2023,Greco2023}.

Continuum flexoelectricity can be framed mathematically as a coupled system of fourth-order partial differential equations, requiring $C^1$-continuous solution fields. Several approaches have been proposed in the literature involving (1) discretization methods based on smooth basis functions, such as the maximum entropy meshless method \cite{Abdollahi2014,Abdollahi2015a,zhuang2020meshfree}, isogeometric approaches \cite{ghasemi2017level, codony2020modeling}, and B-spline based immersed boundary methods \cite{codony2019immersed,codony2021mathematical}, and (2) methods using non-standard finite elements compatible with $C^0$ approximations, such as mixed finite element methods \cite{mao2016mixed,Deng2017,Tian2021}  and interior penalty methods \cite{ventura2020c0}. The formulation and imposition of high-order generalized periodic conditions is not trivial in the general case \cite{mocci2021geometrically}. To impose these conditions in a computational framework, one approach is to encode them in the formulation of a variational principle, for instance using Lagrange multipliers or Nitsche's method \cite{barcelo2022weak,Balcells}.  For open-knot B-Splines and a purely higher-order mechanical problem, these conditions have been imposed by minimizing an algebraic residual function \cite{Schmidt2022}. Here, we present an alternative approach for the solution of high-order boundary value problems of flexoelectricity on a generic unit cell with generalized periodic conditions, based on the construction of a high-order periodic approximation space for the state variables. In this approach, generalized periodicity is built-in the approximation space, and hence the formulation is straightforward and devoid of compatibility conditions on the Lagrange multipliers spaces (LBB) or of penalty parameters. The isogeometric framework is particularly well-suited to the simulation of periodic domains since the construction of a high-order-periodic B-spline basis is trivial, and the cuboidal shape of the fictitious domain can be immediately identified with the unit cell of the architected material. This approach is very elegant and convenient, since it yields an unconstrained boundary value problem in which high-order generalized periodicity conditions are strongly enforced. 

The paper is organized as follows. We first formulate the continuum flexoelectric boundary value problem starting from the  Lifshitz-invariant electromechanical-enthalpy \cite{codony2021mathematical}, with particular attention to the boundary conditions. We then consider an RVE under generalized periodic boundary conditions, and interpret the associated additional degrees of freedom as macroscopic variables. The analogy between the macroscopic response of the RVE and a homogenized material with its macroscopic enthalpy density is discussed. We then describe in detail in Section \ref{Se:Numerical} the proposed numerical approach based on  high-order generalized-periodic approximation spaces for the state variables, as well as the strong enforcement of macroscopic conditions. The proposed method is validated and demonstrated in two examples of functional two-dimensional and three-dimensional flexoelectric metamaterials in Section \ref{Se:examples}.

\section{Standard boundary value problem for flexoelectricity}\label{Preliminary}
\begin{figure}[htb]
    \centering
    \includegraphics[width=0.7\textwidth]{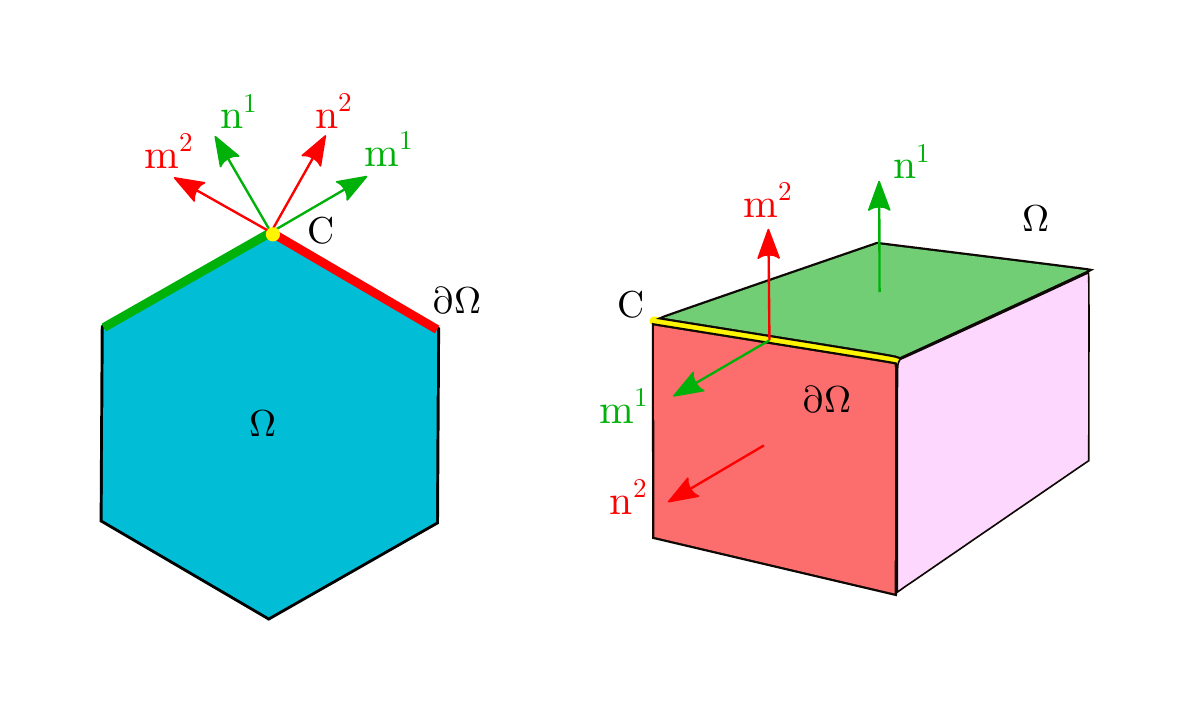}
    \caption{Example of normal and conormal vector for a corner and edge for a 2D case and 3D case respectively.}
    \label{fig:normalconormal}
\end{figure}

Following \cite{codony2021mathematical}, we briefly review the theoretical framework of flexoelectricity used here, considering the displacement field $\Displacement$ and the electric potential $\phi$ as state variables. The total enthalpy $\Pi[\Displacement,\phi]$ of a flexoelectric material occupying a domain  $\Omega$  in $\mathbb{R}^2$ or $\mathbb{R}^3$ is
\begin{equation}\label{eq_Functional}
\Pi[\Displacement,\phi]=
\int_{\Omega}\left(\mathcal{H}^\Omega[\Displacement,\phi]+\mathcal{W}^{\Omega}[\Displacement,\phi]\right)\dd\Omega
+
\int_{\partial\Omega}\mathcal{W}^{\partial\Omega}[\Displacement,\phi]\dd\Gamma
+
\int_{C}\mathcal{W}^C[\Displacement,\phi]\dd\ds,
\end{equation}
where $\partial\Omega$ denotes the domain's boundary, composed by smooth patches joined at sharp edges $C$ (or sharp corners in 2D), as sketched in Fig.~\ref{fig:normalconormal}, The Lifshitz-invariant form of bulk internal enthalpy density in the regime of infinitesimal deformations, accounting explicitly for both direct and converse flexoelectricity, is
\begin{equation}\label{FlexRawEnergy}
\mathcal{H}^\Omega[\Displacement,\phi]=
\frac{1}{2}\strain_{ij}\elast_{ijkl}\strain_{kl}
+\frac{1}{2}\strain_{ij,k}\strGr_{ijklmn}\strain_{lm,n}
-\frac{1}{2}E_{l}\epsilon_{lm}E_{m}
-\frac{1}{2}E_{m,n}M_{mnlk}E_{l,k}
-E_{l}\piezo_{lij}\strain_{ij}
-\frac{1}{2}E_{l}\flexo_{lijk}\strain_{ij,k}
+\frac{1}{2}E_{l,k}\flexo_{lijk}\strain_{ij},
\end{equation}
with the strain $\strain_{ij}(\Displacement)=(u_{i,j}+u_{j,i})/2$ and the electric field $E_l(\phi) = -\phi_{,l}$.
Einstein's summation notation is used, \ie repeated indices sum over the spatial dimensions. Indeces after a comma denote spatial derivatives, i.e. $u_{i,j}=\partial u_i/\partial x_j$. In \eq \eqref{FlexRawEnergy}, $\elast_{ijkl}$ is the elasticity tensor, $\strGr_{ijklmn}$ is the strain-gradient elasticity tensor, $\epsilon_{lm}$ is the dielectricity tensor, $M_{mnlk}$ is the gradient dielectricity tensor, $\piezo_{lij}$ is the piezoelectric tensor and $\flexo_{lijk}$ is the flexoelectric tensor. The components of the material tensors are detailed in Appendix A.

The work from volumetric external loads is
\begin{equation}\label{ExternalBulk}
\mathcal{W}^{\Omega}[\Displacement,\phi]= -b_iu_i+q\phi,
\end{equation}
where $\textbf{b}$ and $q$  represent body force and free charge per unit volume.
External surface loads on the domain boundary $\partial\Omega$ yield the following work per unit area:
\begin{equation}\label{eq_ExternalBoundary}
\mathcal{W}^{\partial\Omega}[\Displacement,\phi]= -t_iu_i-r_i\partial^nu_i+w\phi
+\mathfrak{r}\partial^n\phi,
\end{equation}
where the {traction} $\Traction$, {double traction} $\HighOrderTraction$, {surface charge density} $w$ and {double charge density} $\mathfrak{r}$ are the conjugates of the displacement $\Displacement$, the normal derivative of the displacement $\partial^n(\Displacement)$, the electric potential $\phi$ and the normal derivative of the electric potential $\partial^n\phi$ on $\partial\Omega$, respectively. In addition, {forces per unit length} $\HighOrderJump$ and {electric charges per unit length} $\wp$ arise at the edges $C$ of the boundary, according to the fourth order nature of the associated PDE. Hence, its work per unit length is
\begin{equation}\label{eq_ExternalLine}
\mathcal{W}^C[\Displacement,\phi]= -j_iu_i
+ \wp\phi.
\end{equation}

The boundaries of the domain $\partial\Omega$ and their sharp regions $C$ are split in several sets. Each state variable and its normal derivative determine two disjoint regions: the Dirichlet regions where their value is prescribed, and the Neumann regions where their energy conjugate is prescribed instead. The displacement field determines 
$\partial\Omega=\partial\Omega_u\cup\partial\Omega_t$ and $C=C_u\cup C_j$
, while its normal derivative implies
$\partial\Omega=\partial\Omega_v\cup\partial\Omega_r$.
In turn, the electric potential leads to
$\partial\Omega=\partial\Omega_\phi\cup\partial\Omega_w$ and $C=C_\phi\cup C_\wp$, and its normal derivative determines $\partial\Omega=\partial\Omega_\varphi\cup\partial\Omega_\mathfrak{r}$.
The corresponding boundary conditions are
\begin{subequations}
\begin{align}\label{strgr_classicalBC}
\Displacement - \Displacement^D &= \mathbf{0} \quad\text{on }{\partial\Omega_u},
&&&
\Traction(\Displacement,\phi) - \Traction^N &= \mathbf{0} \quad\text{on }{\partial\Omega_t},
\\\label{strgr_nonlocalBC}
\partial^n(\Displacement) - \bm{v}^D &= \mathbf{0} \quad\text{on }{\partial\Omega_v},
&&&
\HighOrderTraction(\Displacement,\phi) - \HighOrderTraction^N &= \mathbf{0} \quad\text{on }{\partial\Omega_r},
\\\label{elecBC}
\phi-\phi^D &=0\quad\text{on }{\partial\Omega_\phi},
&&&
w(\Displacement,\phi)-w^N &=0\quad\text{on }{\partial\Omega_w},\\
\partial^n(\phi) - \varphi^D&=0 \quad\text{on }{\partial\Omega_\varphi}, 
&&&
\mathfrak{r}(\Displacement,\phi)- \mathfrak{r}^N&=0\quad\text{on }{\partial\Omega_\mathfrak{r}},\\
\Displacement - \Displacement^D &= \mathbf{0}\quad\text{on }{C_u},
&&&
\HighOrderJump(\Displacement,\phi) - \HighOrderJump^N  &= \mathbf{0}\quad\text{on }{C_j},\\
\phi - \phi^D &= 0\quad\text{on }{C_\phi},
&&&
\wp(\Displacement,\phi) - \wp^N &= 0\quad\text{on }{C_\wp},
\end{align}
\end{subequations}
with the prescribed values for the displacement $\Displacement^D$, its normal derivative $\bm{v}^D$, the electric potential $\phi^D$ and its normal derivative $\varphi^D$ at the Dirichlet boundaries, and the prescribed values for the traction
$\Traction^N$, the double traction $\HighOrderTraction^N$, the surface charge density $w^N$, the double surface charge density $\mathfrak{r}^N$, the force per unit length $\HighOrderJump^N$ and the electric charge per unit length $\wp^N$ at the Neumann boundaries. The expressions for $\Traction(\Displacement,\phi)$, $\HighOrderTraction(\Displacement,\phi)$, $w(\Displacement,\phi)$, $\mathfrak{r}(\Displacement,\phi)$, $\HighOrderJump(\Displacement,\phi)$ and $\wp(\Displacement,\phi)$ are derived from the variational principle
\begin{equation}\label{eq_variationalPrinciple}
\left(\Displacement^\text*,\phi^\text*\right)=\arg\min_{\Displacement\in\mathcal{U}}\max_{\phi\in\mathcal{P}}\Pi[\Displacement,\phi],
\end{equation}
which determines the equilibrium states $(\Displacement^*,\phi^*)$. The state variables
$(\Displacement,\phi)\in\mathcal{U}\otimes\mathcal{P}$,
where
\begin{align}
\mathcal{U} &=\left\{\Displacement\in[\mathcal{H}^2(\Omega)]^3\text{ $|$ }
\Displacement-\Displacement^D=\mathbf{0}\text{ on }\partial\Omega_u,\ 
\Displacement-\Displacement^D=\mathbf{0}\text{ on }C_u
\text{ and }
\partial^n\Displacement-\bm{v}^D=\mathbf{0}\text{ on }\partial\Omega_v\right\},
\\
\mathcal{P} &=\left\{\phi\in\mathcal{H}^2(\Omega)\text{ $|$ }
\phi-\phi^D=0\text{ on }\partial\Omega_\phi,\  \phi-\phi^D=0\text{ on }C_\phi \text{ and } \partial^n\phi-\varphi^D=0 \text{ on } \partial\Omega_\varphi \right\},
\end{align}
fulfilling Dirichlet boundary conditions. The resulting Euler-Lagrange equations are derived in \cite{codony2021mathematical} as
\begin{subequations}\label{eq_EulerLagrange}
\begin{align}
\left(\cauchyStress_{ij}(\Displacement,\phi)-\hyperStress_{ijk,k}(\Displacement,\phi)\right)_{,j}+b_i&=0 \qquad\text{in }\Omega,\\
\hfill\left(\electricDisp_{l}(\Displacement,\phi)-\widetilde{D}_{lk,k}(\Displacement,\phi)\right)_{,l}-q&=0\qquad \text{in }\Omega. 
\end{align}
\end{subequations}
and the following expressions are identified
\begin{subequations}{\label{eq_FlexoForces1}}
\begin{align}
t_i(\Displacement,\phi) &=
\left(\cauchyStress_{ij}
-
\hyperStress_{ijk,k}+\surfaceGradient_l\left(n_l\right)\hyperStress_{ijk}n_k\right)n_j
-\surfaceGradient_j\left(\hyperStress_{ijk}n_k\right)
\hspace{-15em}&\text{ on }\partial\Omega,
\\
r_i(\Displacement,\phi) &= \hyperStress_{ijk}n_jn_k
\hspace{-15em}&\text{ on }\partial\Omega, \\
w(\Displacement,\phi) &= -\left(\electricDisp_l-\widetilde{D}_{lk,k}+\nabla^S_ i(n_ i)\widetilde{D}_{lk}n_k\right)n_l+\nabla^S_l\left(\widetilde{D}_{lk}n_k\right)
\hspace{-15em}&\text{ on }\partial\Omega,\\
\mathfrak{r}(\Displacement,\phi) &= -\widetilde{D}_{jk}n_jn_k
\hspace{-15em}&\text{ on }\partial\Omega, \\
j_i(\Displacement,\phi) &= \jump{\hyperStress_{ijk}(\Displacement,\phi)m_jn_k}
\hspace{-15em}&\text{ on }C, \\
\wp(\Displacement,\phi) &= -\jump{\widetilde{D}_{jk}(\Displacement,\phi)m_jn_k}
\hspace{-15em}&\text{ on }C,
\end{align}
\end{subequations}
where \bm{n} is defined on $\partial\Omega$ as the exterior unit normal vector, and \bm{m} represents the conormal vector on $C$, \ie the vector tangent to $\partial\Omega$ and normal to $C$, see  Fig.~\ref{fig:normalconormal}. The jump operator in \eq \eqref{eq_FlexoForces1} is defined on $C$ as $\jump{A}=A^1+A^2$ where superindices 1 and 2 correspond to the two surfaces adjacent to an edge $C$. The symbol
$\nabla^S_j(~) = \nabla_k(~)\left(\delta_{kj}-n_kn_j\right)$
denotes the surface-divergence operator defined on $\partial\Omega$.
In \eqs \eqref{eq_EulerLagrange} and \eqref{eq_FlexoForces1}, the {Cauchy stress} $\CauchyStress(\Displacement,\phi)$, the {double stress} $\HyperStress(\Displacement,\phi)$, the {local electric displacement} $\ElectricDisp(\Displacement,\phi)$ and the {double electric displacement} $\widetilde{\bm{D}}(\Displacement,\phi)$ are defined as the conjugates to the strain $\Strain(\Displacement)$, the strain gradient $\gradient\Strain(\Displacement)$, the electric field $\E(\phi)$ and the electric field gradient $\gradient\E(\phi)$, respectively, as follows:
\begin{subequations}\label{elec_tensors}\begin{align}
\cauchyStress_{ij}(\Displacement,\phi)=\cauchyStress_{ji}(\Displacement,\phi)
&=
\frac{\partial \mathcal{H}^\Omega[\Strain,\gradient\Strain,\E,\gradient\E]}{\partial\strain_{ij}}
=\elast_{ijkl}\strain_{kl}(\Displacement)-\piezo_{lij}E_{l}(\phi)+\frac{1}{2}\flexo_{lijk}E_{l,k}(\phi),
\\
\hyperStress_{ijk}(\Displacement,\phi)=\hyperStress_{jik}(\Displacement,\phi)
&=
\frac{\partial\mathcal{H}^\Omega[\Strain,\gradient\Strain,\E,\gradient\E]}{\partial\strain_{ij,k}}
=\strGr_{ijklmn}\strain_{lm,n}(\Displacement)-\frac{1}{2}\flexo_{lijk}E_{l}(\phi),
\\
\electricDisp_l(\Displacement,\phi)
&=
-\frac{\partial\mathcal{H}^\Omega[\Strain,\gradient\Strain,\E,\gradient\E]}{\partial E_{l}}
=\epsilon_{lm}E_{m}(\phi)+\piezo_{lij}\strain_{ij}(\Displacement)+\frac{1}{2}\flexo_{lijk}\strain_{ij,k}(\Displacement),\\
\widetilde{D}_{kl}(\Displacement,\phi)=
\widetilde{D}_{lk}(\Displacement,\phi)
&=
-\frac{\partial\mathcal{H}^\Omega[\Strain,\gradient\Strain,\E,\gradient\E]}{\partial E_{l,k}}
=M_{mnlk}E_{m,n}(\phi)-\frac{1}{2}\flexo_{lijk}\varepsilon_{ij}(\Displacement).
\end{align}\end{subequations}

The physical stress $\bm{\sigma}$ and the physical electric displacement $\bm{D}$ are deduced from \eq \eqref{eq_EulerLagrange} as

\begin{align}
    \sigma_{ij}&= \cauchyStress_{ij}-\hyperStress_{ijk,k}=
	    {\elast}_{ijkl}\strain_{kl}(\Displacement)-\piezo_{lij}E_{l}(\phi)-\strGr_{ijklmn}\strain_{lm,nk}(\Displacement)+\flexo_{lijk}E_{l,k}(\phi),\nonumber\\
	    D_l &= \electricDisp_l-\widetilde{D}_{lk,k} = 
	    \dielec_{lm}E_{m}(\phi)+\piezo_{lij}\varepsilon_{ij}(\Displacement)-M_{ijlk}E_{i,jk}(\phi)+\flexo_{lijk}\strain_{ij,k}(\Displacement).
\end{align}

\section{Macroscopic conditions for flexoelectric RVE via high-order generalized periodicity}\label{Se:Macroscopic}

In this Section, we state the conditions on the state variables of an RVE that must hold in order to reproduce the bulk response of an infinitely large periodic structure $\Omega^\infty$. The resulting macroscopic state variables and their corresponding macroscopic enthalpy functional are also analyzed.

\begin{figure}[b!]
\centering
    \includegraphics[width=0.8\textwidth]{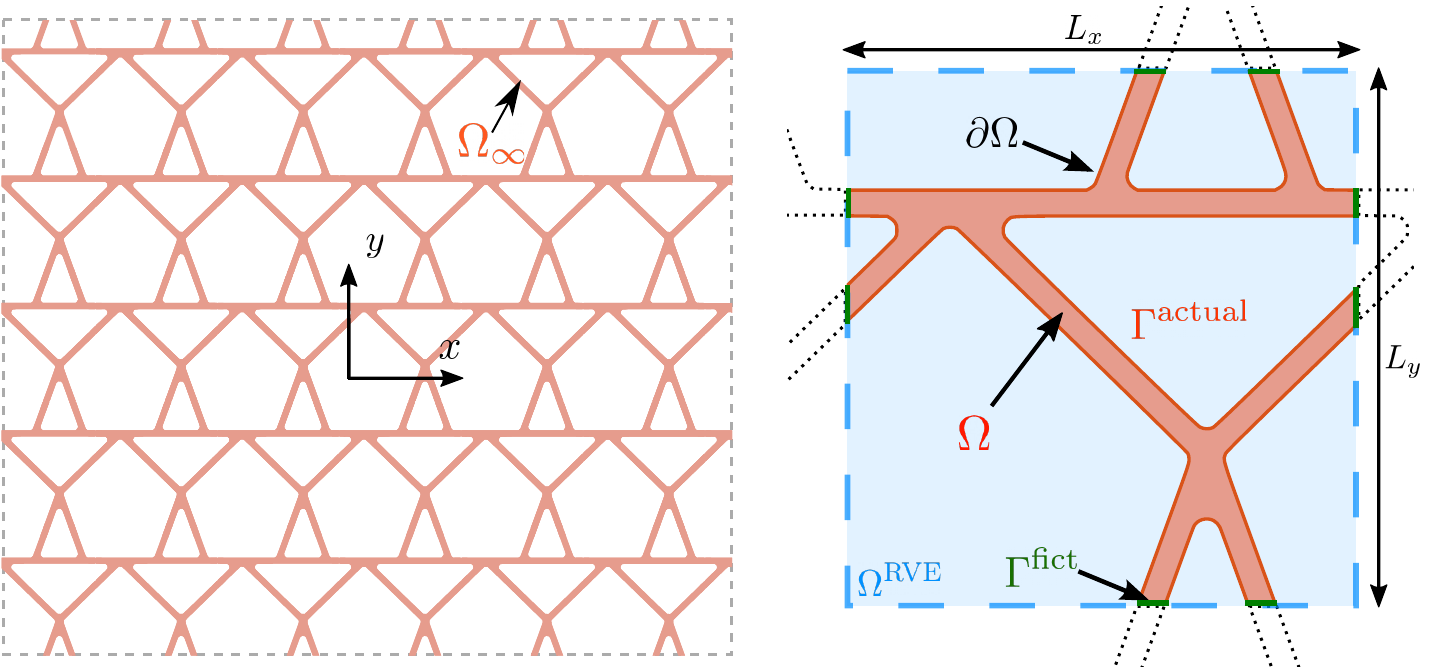}
    \caption[Generic architected material formed by a microstructure $\Omega$ within a unit cell $\Omega^\text{RVE}$]{Generic architected material formed by a microstructure $\Omega$ within a unit cell $\Omega^\text{RVE}$ of dimensions $\left(L_{x},L_{y}\right)$ that is replicated periodically along $x$ and $y$ directions.}
\label{fig_ArchitectedSetup}
\end{figure}

The periodic structure $\Omega^\infty$ is formed by endless concatenation of a unit cell $\Omega^\text{RVE}$ in each spatial dimension. The domain $\Omega$ is the intersection between the periodic structure and the unit cell,  $\Omega=\Omega^\infty\cap \Omega^\text{RVE}$. The boundary of the domain $\partial\Omega$ is split in two parts, $\partial\Omega=\Gamma^{\text{fict}}\cup\Gamma^\text{actual}$ with $\Gamma^\text{fict}=\partial\Omega^\text{RVE}\cap\Omega^\infty$ and $\Gamma^\text{actual}=\partial\Omega^\infty\cap\Omega^\text{RVE}$, see \fig \ref{fig_ArchitectedSetup}. 

For the sake of simplicity, we consider homogeneous Neumann boundary conditions at every physical boundary of the RVE, and disregard volumetric external loads, \ie
\begin{subequations}
\label{eq_AssumptionsHomogenization}
    \begin{align}
        &\partial\Omega_u = \partial\Omega_v = \partial\Omega_\phi = \partial\Omega_\varphi = \emptyset, \\
        &\bm{t}^N=\bm{r}^N=\bm{j}^N=\bm{b}=0,\\
        &w^N=\mathfrak{r}^N=\wp^N=q=0.
    \end{align}
\end{subequations}

\subsection{High-order generalized periodicity conditions}\label{PeriodicityConditions}

\begin{figure}[tb!]
    \centering
    \includegraphics[width=0.95\textwidth]{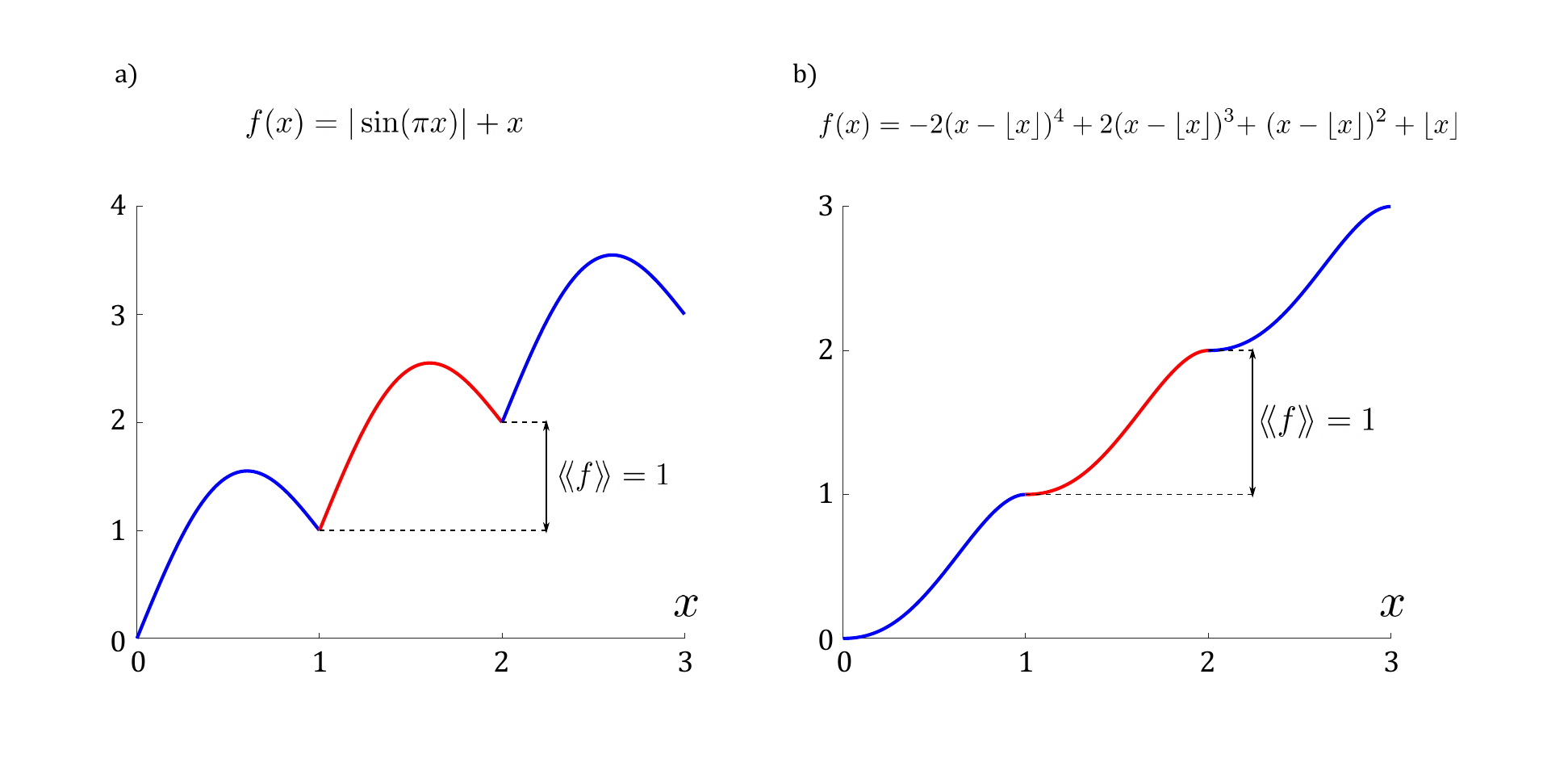}
    \caption{Generic univariate functions $f(x)$ fulfilling (a) low-order ($C^0$) generalized periodicity conditions in \eq \eqref{eq:periodicity1} and (b) high-order ($C^1$) generalized periodicity conditions in \eqs \eqref{eq:periodicity1} and \eqref{eq:periodicity2}.}
    \label{fig:Continuidad1D}
\end{figure}

Generalized periodicity conditions for a generic 1D field $f(x)\in\mathbb{R},x\in\Omega^\text{RVE}=[0,L_x]\subset\mathbb{R}$ are usually stated as
\begin{equation}\label{eq:periodicity1}
    f(L_x)-f(0)=\jumpMacro{f}_x,
\end{equation}
with $\jumpMacro{f}_x\in\mathbb{R}$. Standard periodicity conditions are obtained for $\jumpMacro{f}_x=0$, and \emph{generalized} periodicity conditions otherwise.
In a fourth-order PDE context, this condition is required but insufficient, since the extension of $f$ over $\mathbb{R}$ is required to belong to $H^2(\mathbb{R})$ (\ie it must be at least $C^1$-continuous), which is not necessarily true at $x=m_xL_x$, $m_x\in\mathbb{Z}$. An extra necessary condition is then
\begin{equation}\label{eq:periodicity2}
    \frac{\partial f(L_x)}{\partial x}-\frac{\partial f(0)}{\partial x}=0.
\end{equation}
The difference between low-order periodicity conditions and high-order periodicity conditions is illustrated in \fig \ref{fig:Continuidad1D}. 

The extension of high-order generalized periodicity conditions to higher dimensions is trivial. Let us consider a cuboidal unit cell $\Omega^\text{RVE}=[0,L_x]\times[0,L_y]\times[0,L_z]\in\mathbb{R}^3$ as depicted in \fig \ref{fig_ArchitectedSetup} for the 2D case. The high-order generalized periodicity conditions of the mechanical and electrical fields are:
\begin{subequations}\label{GeneralizedPeriodicityConditions}
\begin{align}
&\Displacement(\zeta=L_\zeta)-\Displacement(\zeta=0)=\jumpMacro{\Displacement}_\zeta, &&&
&\phi(\zeta=L_\zeta)-\phi(\zeta=0)=\jumpMacro{\phi}_\zeta,\\
\label{eq:High-orderContinuity}&\frac{\partial\Displacement(\zeta=L_\zeta)}{\partial \zeta}-\frac{\partial\Displacement(\zeta=0)}{\partial\zeta}=0, &&&
&\frac{\partial\phi(\zeta=L_\zeta)}{\partial\zeta}-\frac{\partial\phi(\zeta=0)}{\partial\zeta}=0,
\end{align}
\end{subequations}
for $\zeta=\{x,y,z\}$. 
\subsection{High-order equilibrium conditions}
\label{Se:High-order equilibrium conditions}
In addition to the continuity conditions stated above, we need to impose high-order equilibrium of the solution fields across the RVE boundaries, that is \cite{barcelo2022weak},
\begin{align}
\label{eq:macrotraction}
    \bm{t}(\zeta=L_\zeta)+\bm{t}(\zeta=0)&=0, &&&
    \bm{r}(\zeta=L_\zeta)-\bm{r}(\zeta=0)&=0, &&&
    \bm{j}(\zeta=L_\zeta)+\bm{j}(\zeta=0)&=0, \nonumber \\
    w(\zeta=L_\zeta)+w(\zeta=0)&=0, &&&
    \mathfrak{r}(\zeta=L_\zeta)-\mathfrak{r}(\zeta=0)&=0, &&& 
    \wp(\zeta=L_\zeta)+\wp(\zeta=0)&=0.
\end{align}

Note that \eq \eqref{eq:macrotraction} is actually required even if \eqs \eqref{eq:High-orderContinuity} hold, since the Neumann quantities $(\bm{t},\bm{r},\bm{j},w,\mathfrak{r},\wp)$, depend also on second and third-order derivatives of the state variables $(\Displacement,\phi)$ that are not periodic in general.


\subsection{Macroscopic kinematics}
\label{Se:MacroscopicKinematics}
The additional unknowns $\jumpMacro{f}_{\zeta }$ arising from the generalized periodicity conditions represent the jump (or difference between boundaries) on $\Omega^\text{RVE}$ of the field $f(\x)$ with $\x\in\Omega$ along the Cartesian direction $\zeta$. From a macroscopic point of view, they are regarded as the state variables that characterize the macroscopic behavior of an homogenized field on $\Omega^\text{RVE}$, regardless of the evolution of $f(\x)$ within $\Omega$ or even the shape of $\Omega$. 

We are interested in capturing the macroscopic behaviour of the state variables: displacement $\Displacement(\x)$ through the macroscopic displacement gradient $\overline{\gradient \Displacement}$ with nine unknowns in $\mathbb{R}^3$ and electric potential $\phi(\x)$ through the macroscopic electric field $\overline{\gradient\E}$ with three unknowns. 
By taking into account these macroscopic entities and the dimensions of $\Omega^\text{RVE}$, the displacement $\Displacement$ and electric potential $\phi$ are split into microscopic  and macroscopic contributions,
\begin{subequations}
\label{eq:StateDecomposition}
    \begin{align}
        &\Displacement(\bm{x})=\Displacement^\text{P}(\bm{x})+ \overline{\nabla\Displacement}\cdot \bm{x}, \\
        &\phi(\bm{x})=\phi^\text{P}(\bm{x})- \overline{\nabla\E}\cdot \bm{x}, 
    \end{align}
\end{subequations}
where $\Displacement^\text{P}(\bm{x})$ and $\phi^\text{P}(\bm{x})$ are periodic functions that fulfill

\begin{subequations}\label{GeneralizedPeriodicityConditionsPeriodicSubspace}
\begin{align}
\label{eq:DisplacementPeriodicity}
\Displacement^\text{P}(\zeta=L_\zeta)&-\Displacement^\text{P}(\zeta=0)=0, &&&
\frac{\partial\Displacement^\text{P}(\zeta=L_\zeta)}{\partial \zeta}&-\frac{\partial\Displacement^\text{P}(\zeta=0)}{\partial\zeta}=0,\\
\label{eq:PotentialPeriodicity}\phi^\text{P}(\zeta=L_\zeta)&-\phi^\text{P}(\zeta=0)=0, &&&
\frac{\partial\phi^\text{P}(\zeta=L_\zeta)}{\partial\zeta}&-\frac{\partial\phi^\text{P}(\zeta=0)}{\partial\zeta}=0,
\end{align}
\end{subequations}
and the macroscopic displacement gradient and electric field are
\begin{equation}
    \overline{\nabla\Displacement}=\left[
    \begin{array}{ccc}
         \jumpMacro{u_x}_x/L_{x} & 
         \langle\!\langle u_y\rangle \!\rangle _x/L_{x} & 
         \jumpMacro{u_z}_x/L_{x}  \\
         \jumpMacro{u_x}_y/L_{y} & 
         \langle\!\langle u_y\rangle \!\rangle_y/L_{y} & 
         \jumpMacro{u_z}_y/L_{y}  \\
         \jumpMacro{u_x}_z/L_{z} & 
         \langle\!\langle u_y\rangle \!\rangle_z/L_{z} & 
         \jumpMacro{u_z}_z/L_{z}  \\
    \end{array}\right],
\end{equation}
\begin{equation}\label{eq:macroEF}
    \overline{\E}=-\overline{\nabla\phi}=\left[
    \begin{array}{c}
         -\jumpMacro{\phi}_x/L_{x} \\
         -\jumpMacro{\phi}_y/L_{y} \\
         -\jumpMacro{\phi}_z/L_{z} \\
    \end{array}\right].
\end{equation}
In turn, the macroscopic displacement gradient can be uniquely decomposed into its symmetric and antisymmetric parts as
\begin{subequations}
\label{eq:MacroDisplacementDecomposition}
\begin{align}
&\overline{\nabla\Displacement}=
    \overline{\Strain}+\overline{\textbf{W}},\label{eq:MacroDisplacementDecomposition2}\\
    &\overline{\Strain}=\frac{1}{2}\left(
    \overline{\nabla\Displacement}+\overline{\nabla\Displacement}^T\right),\\
    &\overline{\textbf{W}}=\frac{1}{2}\left(
    \overline{\nabla\Displacement}-\overline{\nabla\Displacement}^T\right).
\end{align}
\end{subequations}
The macroscopic strain $\overline{\Strain}$ and the macroscopic infinitesimal rotation $\overline{\textbf{W}}$ (spin) are constant tensors that characterize the macroscopic (homogenized) kinematics associated to a generalized periodic displacement field $\Displacement(\x)$ defined on $\Omega$. They are invariant with respect to the specific choice of RVE given a periodic material, e.g.~number of repetitions of the unit cell or translations of the RVE boundary on the material frame. The deformation is captured by  $\overline{\Strain}$, whereas $\overline{\textbf{W}}$ represents a rigid-body spin that leaves the energy functional unchanged. Computationally, this spin needs to be fixed to have a unique solution. Without loss of generality, we consider $\overline{\textbf{W}}=\bm{0}$, that is, we assume symmetric macroscopic displacement gradients $(\overline{\nabla\Displacement}=
\overline{\Strain})$.

\begin{remark}
Note that assuming $\Displacement,\phi$ to be generalized-periodic implies that the macroscopic state variables $\overline{\Strain},\overline{\E}$ are constant tensors. Therefore, the macroscopic strain gradient and electric field gradient vanish. Loading cases on a periodically-arranged architected material yielding non-vanishing macroscopic strain gradient and electric field gradients are out of the scope of this work.
\end{remark}

\begin{remark}
Equation \eqref{eq:StateDecomposition} implies that the strain field and electric field are decomposed as 
\begin{subequations}
\label{eq:StrainElectricDecomposition}
    \begin{align}
        &\Strain = \Strain \left(\Displacement^{\text{P}}\right) + \overline{\Strain} \\
        &\E = \E\left(\phi^{\text{P}}\right) + \overline{\E}.
    \end{align}
\end{subequations}
\end{remark}

\subsection{Macroscopic enthalpy functional and conjugate variables}
\label{Se:energyequivalence}
Since the macroscopic response of the architected structure is uniquely characterized by the macroscopic state variables, our goal is to rationalize the existence of a macroscopic enthalpy functional depending on macroscopic state variables only. Such functional should fulfill the condition that the variation of the actual bulk enthalpy of the system or microscopic bulk enthalpy ($\Pi^\text{b}$) in one unit cell is equivalent to the variation of the enthalpy of an homogeneous media of size $L_x$, $L_y$ and $L_z$ or macroscopic enthalpy ($\overline{\Pi}$):
\begin{equation}\label{Balance}
\delta\Pi^\text{b}[\Strain,\nabla\Strain,\E,\nabla\E]=
\delta\overline{\Pi}\left[\overline{\Strain},\overline{\E}\right].
\end{equation}

In order to find the functional described in \eq \eqref{Balance}, let us consider the high-order equilibrium conditions from \eq \eqref{eq:macrotraction} in weak form as
\begin{equation}
\label{eq:BoundaryFict}
    0=\int_{\Gamma_{\text{fict}}}\left(-t_i\delta u_i^\text{P}-r_i\partial^n\delta u_i^\text{P}+w\delta\phi^\text{P}+\mathfrak{r}\partial^n\delta \phi^\text{P}\right)\dd\Gamma + \int_{C_{\text{fict}}}\left(-j_i\delta u_i^\text{P}+\wp \delta\phi^\text{P}\right)\dd l,
\end{equation}
for all admissible periodic test functions whose functional space will be defined later in Section \ref{Se:BoundaryValueProblem}. 
Assuming homogeneous microscopic Neumann conditions on the physical boundaries of the RVE, c.f.~\eq \eqref{eq_AssumptionsHomogenization}, \eq \eqref{eq:BoundaryFict} is extended to $\Gamma_{\text{fict}}\cup\Gamma_{\text{act}}=\partial\Omega$, $C_{\text{fict}}\cup C_{\text{act}}=C$, so
\begin{equation}
\label{eq:BoundaryFictAct}
    0=\int_{\partial \Omega}\left(-t_i\delta u_i^\text{P}-r_i\partial^n\delta u_i^\text{P}+w\delta\phi^\text{P}+\mathfrak{r}\partial^n\delta \phi^\text{P}\right)\dd\Gamma + \int_{C}\left(-j_i\delta u_i^\text{P}+\wp \delta\phi^\text{P}\right)\dd l.
\end{equation}
Upon integration by parts, and invoking the divergence theorem and the surface divergence theorem, as done in \cite{codony2021mathematical} and considering the strong form of the problem in \eq \eqref{eq_EulerLagrange} with zero source terms, c.f.~\eqref{eq_AssumptionsHomogenization}, we obtain
\begin{equation}
    0=\int_{\Omega}\left(\cauchyStress_{ij}\delta\varepsilon_{ij}\left(\Displacement^{\text{P}}\right)+\hyperStress_{ijk}\delta\varepsilon_{ij,k}\left(\Displacement^{\text{P}}\right)-\widehat{D}_l\delta E_l\left(\phi^{\text{P}}\right) - \widetilde{D}_{lm}\delta E_{l,m}\left(\phi^{\text{P}}\right)\right)\dd\Omega.
\end{equation}
By using the decomposition of the strain and electric field in \eq \eqref{eq:StrainElectricDecomposition} we have

\begin{align}
\label{eq:LastVariationalEquation}
    &\int_{\Omega}\left(\cauchyStress_{ij}\delta\varepsilon_{ij}+\hyperStress_{ijk}\delta\varepsilon_{ij,k}-\widehat{D}_l\delta E_l - \widetilde{D}_{lm}\delta E_{l,m}\right)\dd\Omega=
    \int_{\Omega}\cauchyStress_{ij}\delta\overline{\varepsilon}_{ij}\dd\Omega-\int_{\Omega}\hat{D}_l\delta\overline{E}_l\dd\Omega.
\end{align}
Comparing \eq \eqref{Balance} with \eq \eqref{eq:LastVariationalEquation}, the variation of microscopic bulk enthalpy corresponds to the integral over $\Omega$ of the bulk internal enthalpy density variation in \eq \eqref{FlexRawEnergy} as
\begin{equation}
    \delta\Pi^\text{b}[\Strain,\nabla\Strain,\E,\nabla\E]=
    \int_{\Omega}\left(\delta\mathcal{H}^\Omega[\Displacement,\phi]\right)\dd\Omega=
    \int_\Omega\left(\cauchyStress_{ij}\delta\strain_{ij}-\electricDisp_l\delta E_l+\hyperStress_{ijk}\delta\strain_{ij,k}-
    \widetilde{D}_{lm}\delta E_{l,m}\right)\dd \Omega,
\end{equation}
and the variation of the macroscopic enthalpy is
\begin{equation}\label{eq:VariationPiMacro}
    \delta \overline{\Pi}[\overline{\Strain},\overline{\E}] = \int_{\Omega^{\textnormal{RVE}}}\left(\widehat{\sigma}_{ij}\delta \overline{\strain}_{ij} - \widehat{D}_l\delta\overline{E}_l\right)\dd \Omega=|\Omega^{\textnormal{RVE}}|\left(\overline{\sigma}_{ij}\delta \overline{\strain}_{ij} - \overline{D}_l\delta\overline{E}_l\right),
\end{equation}
where the term $\overline{\bm{\sigma}}$ in \eq \eqref{eq:LastVariationalEquation} corresponds to the macroscopic stress and $\overline{\bm{D}}$ is the macroscopic electric displacement. The second equality in \eq \eqref{eq:VariationPiMacro} holds by considering that the macroscopic quantities do not depend on the position $\bm{x}$, and $|\Omega^{\textnormal{RVE}}|=L_xL_yL_z$ corresponds to the macroscopic volume of the RVE. The macroscopic stress and macroscopic electric displacement are conjugates of the macroscopic strain and macroscopic electric field, respectively, and they are defined as

\begin{align}\label{eq:DefMacrostress2}
    &\overline{\sigma}_{ij}=\frac{1}{|\Omega^{\textnormal{RVE}}|}\int_{\Omega}\cauchyStress_{ij}\dd\Omega, 
    &\overline{D}_{l}=\frac{1}{|\Omega^{\textnormal{RVE}}|}\int_{\Omega}\electricDisp_l\dd\Omega.
\end{align}
As a result, the macroscopic stress is nothing but the macroscopic average of the microscopic Cauchy stress over $\Omega$, and the macroscopic electric displacement is the macroscopic average of the microscopic local electric displacement over $\Omega$. \eqs \eqref{eq:DefMacrostress2} can be regarded as an extension of the Hill-Mandel theorem \cite{hill1963elastic,hill1967essential} to high-order electromechanics.

The weak equation that generalized periodic state variables $\Displacement(\bm{x})$, $\phi(\bm{x})$ must fulfill in order to reproduce the electromechanical state of an infinitely large periodic structure over $\Omega^\infty$ is
\begin{equation}
\label{eq:FinalWeakForm}
    \int_\Omega\left(\cauchyStress_{ij}\delta\strain_{ij}-\electricDisp_l\delta E_l+\hyperStress_{ijk}\delta\strain_{ij,k}-
    \widetilde{D}_{lm}\delta E_{l,m}\right)\dd \Omega =
    |\Omega^{\textnormal{RVE}}|\left(\overline{\sigma}_{ij}\delta \overline{\strain}_{ij} - \overline{D}_l\delta\overline{E}_l\right).
\end{equation}

\begin{remark}
From \eq \eqref{eq:DefMacrostress2}, we see that the macroscopic stress has the same symmetry as the microscopic Cauchy stress. This is in agreement with the aforementioned fact that the macroscopic enthalpy functional must not depend on the macroscopic spin $\overline{\textbf{W}}$. Indeed, if $\delta\overline{\Strain}$ was replaced by $\delta\overline{\nabla\Displacement}$ in \eq \eqref{eq:VariationPiMacro}, after decomposing it via \eq\eqref{eq:MacroDisplacementDecomposition2} we would obtain the same macroscopic enthalpy plus an additional term $|\Omega^{\textnormal{RVE}}|\overline{\sigma}_{ij}\delta\overline{\textnormal{W}}_{ij}$. It is easy to see that this product always vanishes, regardless of the value of $\delta\overline{\textbf{W}}$, since $\overline{\bm{\sigma}}$ is symmetric and $\delta\overline{\textbf{W}}$ is antisymmetric.
\end{remark}

\begin{remark}
The expressions of the macroscopic stress and electric displacement derived in \eq\eqref{eq:DefMacrostress2} are different from those derived in Ref.~\cite{Balcells}. Both expressions are equivalent as reported in Appendix B, but our definition here does not depend on quantities defined over fictitious RVE boundaries.
\end{remark}

\subsection{Boundary value problem for flexoelectric RVE}
\label{Se:BoundaryValueProblem}
In the generalized periodicity framework described in above Sections, $\overline{\Strain}$ and $\overline{\E}$ are additional state variables, whose components can be specified a priori, or obtained as result of the boundary value problem. To this end, their components are split in two disjoint sets, the macroscopic Dirichlet components

\begin{subequations}\label{eq:MacroscopicDirichlet}
\begin{align} \label{eq:MacroDirichletStrain}
&\overline{\varepsilon}_{ij}=\overline{\varepsilon}^D_{ij}&&& &\text{for }(i,j)\in\mathcal{I}^\Strain, \\
\label{eq:MacroDirichletElectricField}&\overline{E}_l=\overline{E}_l^D &&& &\text{for }l\in\mathcal{I}^\E,
\end{align}
\end{subequations}
and the macroscopic Neumann components \begin{subequations}\label{eq:MacroscopicNeumann}
\begin{align}
&\overline{\sigma}_{ij}=\overline{\sigma}^N_{ij}&&& &\text{for }(i,j)\in\{1,2,3\}\times\{1,2,3\}\setminus\mathcal{I}^\Strain, \\
&\overline{D}_{l}=\overline{D}^N_{l}&&& &\text{for }l\in\{1,2,3\}\setminus\mathcal{I}^\E,
\end{align}
\end{subequations}
where $\mathcal{I}^\Strain\subseteq\{1,2,3\}\times\{1,2,3\}$ such that if $(i,j)\in \mathcal{I}^\Strain$ then $(j,i)\in \mathcal{I}^\Strain$, and $\mathcal{I}^\E\subseteq\{1,2,3\}$, are the subsets of components where macroscopic Dirichlet conditions are applied. 
The right hand side of \eq \eqref{eq:FinalWeakForm} is split in Dirichlet and Neumann components accordingly. Macroscopic Dirichlet conditions are applied strongly, projecting the solution to a functional space that satisfies macroscopic Dirichlet conditions. On the other hand, macroscopic Neumann conditions remain in the weak form. The macroscopic Neumann conditions can be seen as the \emph{natural} macroscopic conditions of the boundary value problem, since the neglect of the macroscopic Neumann term leads to homogeneous macroscopic Neumann conditions.

The weak form of the flexoelectric generalized problem is

\begin{equation}
\emph{Find $\left(\Displacement^\textnormal{P},\phi^\textnormal{P},\overline{\Strain},\overline{\E}\right)\in\mathcal{U}^\textnormal{P}\otimes\mathcal{P}^\textnormal{P}\otimes\overline{\mathcal{U}}^\textnormal{D}\otimes\overline{\mathcal{P}}^\textnormal{D}$}\emph{ such that}
\nonumber
\end{equation}
\begin{equation}
\int_{\Omega}\left(\delta\varepsilon_{ij}\widehat{\sigma}_{ij} +
\delta\varepsilon_{ij,k}\widetilde{\sigma}_{ijk}
-\delta E_l\widehat{D_l}-\delta E_{l,m}\widetilde{D}_{lm}\right)\dd\Omega
=|\Omega^\text{RVE}|\overline{\sigma}^N_{ij}\delta\overline{\varepsilon}_{ij}-|\Omega^\text{RVE}|\overline{D}^N_{l}\delta\overline{E}_{l}
,\nonumber
\end{equation}
\begin{equation}
\label{eq:FinalWeakFormGeneralized}
\forall\left(\delta\Displacement^\textnormal{P},\delta\phi^\textnormal{P},\delta\overline{\Strain},\delta\overline{\E}\right)\in\mathcal{U}^\textnormal{P}\otimes\mathcal{P}^\textnormal{P}\otimes\overline{\mathcal{U}}^0\otimes\overline{\mathcal{P}}^0,
\end{equation}
where
\begin{subequations}
\begin{align}
\mathcal{U}^\textnormal{P} &=\left\{\Displacement^\text{P}\in[\mathcal{H}^2(\Omega)]^3\text{ $|$ }
\text{Eq. } \eqref{eq:DisplacementPeriodicity} \text{ holds (high-order periodicity on $\Displacement^\text{P}$)}
\right\},
\\
\mathcal{P}^\textnormal{P} &=\left\{\phi^\text{P}\in\mathcal{H}^2(\Omega)\text{ $|$ }
\text{Eq. } \eqref{eq:PotentialPeriodicity} \text{ holds (high-order periodicity on $\phi^\text{P}$)}
\right\}, \\
\overline{\mathcal{U}}^\textnormal{D} &=\left\{\overline{\Strain}\in\left[\mathbb{R}^3\right]^2\text{ $|$ }
\text{Eq. } \eqref{eq:MacroDirichletStrain} \text{ holds (macroscopic Dirichlet conditions on $\overline{\Strain}$) and }
\overline{\varepsilon}_{ij}=\overline{\varepsilon}_{ji}
\right\}, \\
\overline{\mathcal{P}}^\textnormal{D} &=\left\{\overline{\E}\in\mathbb{R}^3\text{ $|$ }
\text{Eq. } \eqref{eq:MacroDirichletElectricField} \text{ holds (macroscopic Dirichlet conditions on $\overline{\E}$)}
\right\}, \\
\overline{\mathcal{U}}^0 &=\left\{\delta\overline{\Strain}\in\left[\mathbb{R}^3\right]^2\text{ $|$ }
\delta\overline{\varepsilon}_{ij}=0\text{ for }(i,j)\in\mathcal{I}^\Strain
\right\}, \\
\overline{\mathcal{P}}^0 &=\left\{\delta\overline{\E}\in\mathbb{R}^3\text{ $|$ }
\delta \overline{E}_{i}=0\text{ for }i\in\mathcal{I}^\E
\right\}.
\end{align}
\end{subequations}
\section{Numerical approach: High-order generalized-periodic approximation spaces}\label{Se:Numerical}

In this section we construct a high-order generalized-periodic approximation space by modifying in an elegant way the immersed boundary B-spline-based approach described in \cite{codony2019immersed}. Section \ref{Se:ImmersedBoundaryApproach} summarizes the framework in \cite{codony2019immersed} which considers a functional space with high-order continuity. In Sections \ref{Se:PeriodicBasis} and \ref{Se:GeneralizedPeriodicBasis} we tailor the aforementioned space for generalized-periodic functions. Finally Section \ref{Se:EnforcementMacroscopicKinematics} restrict the approximation space to strongly enforce the macroscopic conditions defined in Section \ref{Se:BoundaryValueProblem}.

\subsection{High-order approximation space:  Immersed boundary B-spline approach}
\label{Se:ImmersedBoundaryApproach}
The discretization of the weak form in \eq  \eqref{eq:FinalWeakFormGeneralized} requires high-order generalized-periodic approximation spaces for the displacement and electric potential in $[\mathcal{H}^2(\Omega)]^3$ and $\mathcal{H}^2(\Omega)$ respectively. Here, following \cite{codony2019immersed}, we consider the immersed boundary B-spline approach \cite{deBoor2001,codony2021mathematical}. Let us consider B-spline basis functions, that is, piece-wise polynomial functions with $C^{q-1}$ continuity, being $q\geq 2$ the degree of approximation. The  uniform univariate B-spline basis $\{B_i^q\}_{i=0}^{n_\xi-1}$ is defined in a parametric space $\xi\in[0,n_\xi+q]$ with the following recursive formula:

\begin{equation}\begin{aligned}\label{eq_BSpline}
B_i^0(\xi)=
\begin{cases}
1 & \xi_i\leq\xi<\xi_{i+1}\\
0 & \text{otherwise}
\end{cases};
&\quad&
B_i^k(\xi)=\frac{\xi-\xi_i}{\xi_{i+k}-\xi_i}B_i^{k-1}(\xi)+\frac{\xi_{i+k+1}-\xi}{\xi_{i+k+1}-\xi_{i+1}}B_{i+1}^{k-1}(\xi);
&\quad&
\begin{aligned}
&k=1,\dots,q\\
&i= 0,\dots,n_\xi+q-k-1,
\end{aligned}
\end{aligned}
\end{equation}
where $\{\xi_i\}^{n_\xi-1}_{i=0}$ are the so-called knot points, here assumed to be not repeated and equispaced (see \fig \ref{fig_PeriodicBasis}). B-splines are defined in a multivariate space by the tensor product of univariate ones, i.e.,

\begin{equation}\label{eq_trivariate}
	B_\toVect{i}^q([\xi,\eta,\tau])
	=
	B_{i_{\xi}}^q(\xi) B_{i_{\eta}}^q(\eta) B_{i_{\tau}}^q(\tau) ;\quad
	i_{\xi}  = 0,\dots,n_{\xi}-1; \quad
	i_{\eta} = 0,\dots,n_{\eta}-1; \quad
	i_{\tau} = 0,\dots,n_{\tau}-1.
\end{equation}

From now on, we omit the superscript $q$ for convenience. We consider a uniform Cartesian mesh $\Omega_\square$ embedding the domain $\Omega$ in an unfitted way, with elements of size $(h_x, h_y, h_z)$.
In the physical space, the unkowns of the problem $\Displacement$ and $\phi$ are approximated as \cite{codony2019immersed,codony2021mathematical}
\begin{subequations}\label{eq:PhysicalSpace}
    \begin{align}
        &[\Displacement(\bm{x})]_d\approx [\Displacement^h(\bm{x})]_d=N_i(\bm{x})a^u_{id}=B_i(\xi,\eta,\tau)a^u_{id},\\
        &\phi(\bm{x})\approx\phi^h(\bm{x})=N_i(\bm{x})a^\phi_{i}=B_i(\xi,\eta,\tau)a^\phi_{i},
    \end{align}
\end{subequations}
where $N=[B\circ \varphi^{-1}]$, $\varphi(\xi,\eta,\tau)= 
[h_x \xi; h_y \eta; h_z \tau]$ is the geometrical map which maps each point in the parametric space to a given point in the physical one and $\{\bm{a}^u,\bm{a}^\phi\}$ are the degrees of freedom of $u^h$ and $\phi^h$.

\begin{figure}[h!]
    \centering
    \includegraphics[width=0.8\textwidth]{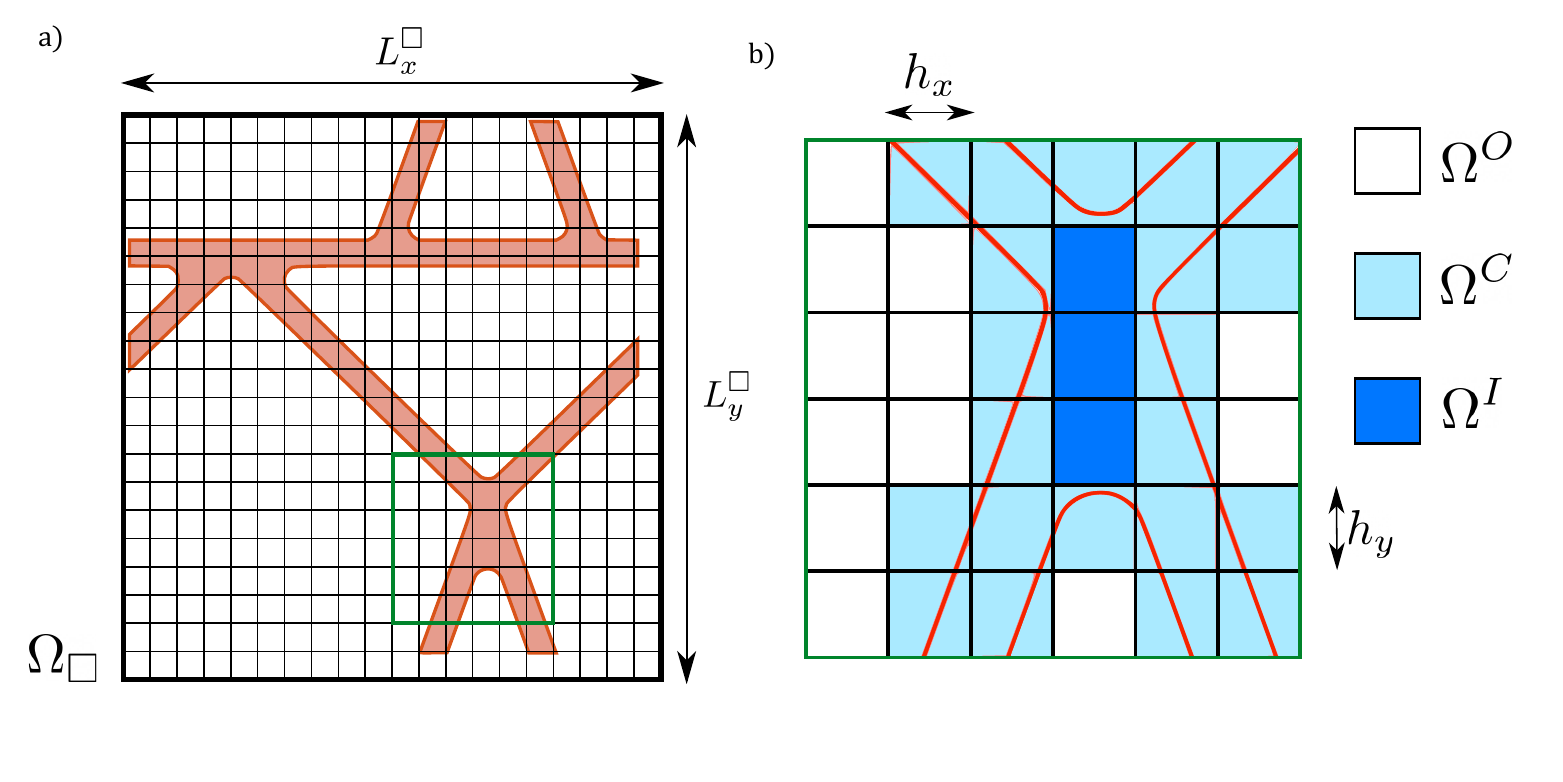}
    \caption{a) Embedded mesh $\Omega_\square$ of size $L_x^\square$ and $L_y^\square$. The number of elements per dimension is 20: $n_x=n_y=19$ and the shift is: $s_x=0.9$ and $s_y=0.1$. b) Zoom in of the mesh, outer cells are depicted in white, cut cells in light blue and inner cells in blue.}
    \label{fig:malla}
\end{figure}

Each cell of the mesh is classified in one of three disjoint groups: inner cells $\Omega^I$, completely contained in the domain $\Omega$, outer cells $\Omega^O$, with no intersection with $\Omega$, and cut cells $\Omega^C$, that are cells intersected by the boundary $\partial\Omega$ (see \fig \ref{fig:malla}).
Cut cells require special treatments regarding numerical integration and cut-cell stabilization, as common in all unfitted methods. Further details about the numerical approach can be found in \cite{codony2019immersed,codony2021mathematical}.
\subsection{High-order periodic space: Periodic basis}
\label{Se:PeriodicBasis}

Let us consider an architected material with unit cell $\Omega^\text{RVE}=[0,L_{x}]\otimes[0,L_{y}]\otimes[0,L_{z}]\subset\mathbb{R}^3$ (see \fig \ref{fig_ArchitectedSetup}).
In order to accommodate a periodic approximation space, the Cartesian mesh $\Omega_\square$ must have element sizes $(h_x, h_y, h_z)$ fulfilling $L_\zeta/h_\zeta = n_\zeta\in \mathbb{N}^+, \zeta=\{x,y,z\}$, where $n_\zeta$ is the number of cells along the $\zeta$-th dimension.
The size of $\Omega_\square$ is denoted by $(L_x^\square, L_y^\square, L_z^\square)$, with
\begin{equation}
    L_\zeta^\square =
    \left\{
    \begin{array}{lc}
    h_\zeta \cdot n_\zeta & \text{  if } L_\zeta^\square = L_\zeta, \\
    h_\zeta \cdot (n_\zeta+1)& \text{  if } L_\zeta^\square > L_\zeta.
    \end{array}
    \right.
\end{equation}
Note that if the former case holds for every dimension, then $\Omega_\square \equiv \Omega^\text{RVE}$.
However, from now on, we restrict ourselves to the more general case where the latter holds for every dimension, yielding $\Omega_\square \supset \Omega^\text{RVE}$. The resulting embedding domain $\Omega_\square$ spans $[ h_\zeta \cdot (s_\zeta-1) , L_\zeta + h_\zeta \cdot s_\zeta ]$  along the $\zeta$-th dimension, with $s_\zeta \in (0,1)$ being an arbitrary shift parameter for each dimension (see \fig\ref{fig:malla}).

\begin{figure}[t!]
\centering
    \includegraphics[width=0.6\textwidth]{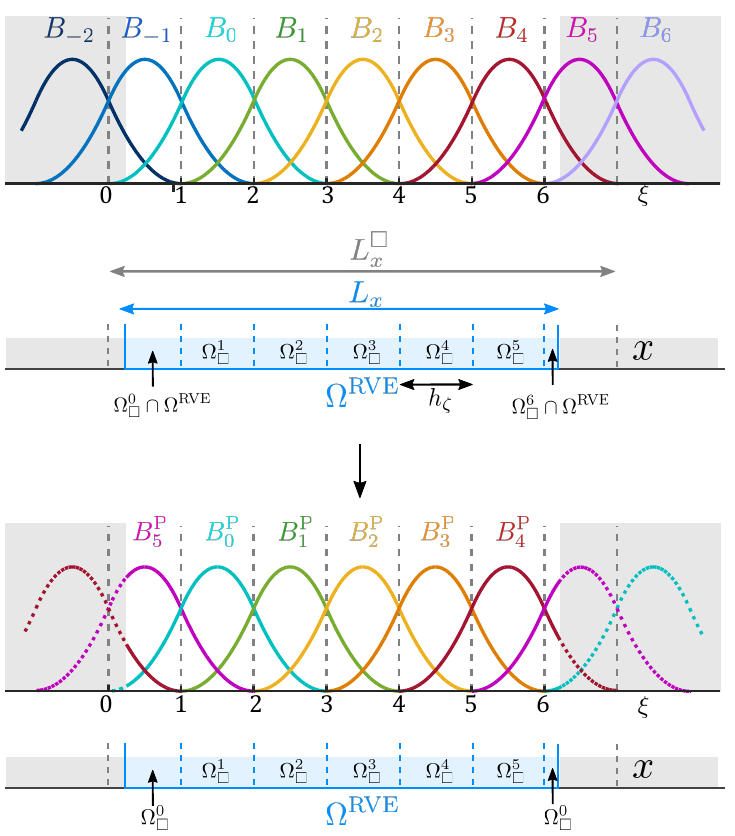}
    \caption[Univariate periodic basis of degree $q=2$]
    {Univariate periodic basis of degree $p=2$. Originally (top), the functional space is spanned by nine B-spline bases $B_i(\xi),i=-2,\dots,6$, defined onto a mesh of $L_x/h_x=n_x=6$ cells. Note that $\Omega^\text{RVE}$ does not coincide with $\Omega_\square$, and therefore cut cells ($\Omega_0^\square\cap\Omega^\text{RVE}$ and $\Omega_6^\square\cap\Omega^\text{RVE}$) are generated. In order to create a high-order-periodic functional space (bottom), the basis functions at a distance $n_x$ on $\Omega^\text{RVE}$ are identified with the same degree of freedom, yielding a functional space spanned by only six periodic B-Spline bases $B_i^\text{P}
   (\xi),i=0\dots5$. The periodic nature of the basis implies periodicity on $\Omega^\text{RVE}$ too, which does not have cells cut by the \emph{periodic boundary} $\partial\Omega^\text{RVE}$ anymore.}
\label{fig_PeriodicBasis}
\end{figure}

Once the mesh is suitably defined, the periodic B-spline basis ${B}^\text{P}_\toVect{i}(\xi,\eta,\tau)$ is obtained from a uniform B-spline basis $B_\toVect{i}(\xi,\eta,\tau)$ by identifying each basis function and its corresponding periodic images at a Cartesian distance of
$(m_{x}n_{x}, m_{y}n_{y}, m_{z}n_{z})$, $m_{x},m_{y},m_{z}\in\mathbb{Z}$
with the same degree of freedom, as illustrated in \fig\ref{fig_PeriodicBasis} for the univariate case.
The periodic nature of $B_\toVect{i}^\text{P}(\xi,\eta,\tau)$ implies high-order periodicity on $\Omega^\text{RVE}$, which is now discretized by uncut cells $\Omega_\square^c$ as shown in \fig \ref{fig_PeriodicBasis}. 
The architected structure of $\Omega\subseteq\Omega^\text{RVE}$ is immersed into the periodic mesh, generating an unfitted discretization with cut cells intersected by $\partial\Omega$, but not by $\partial\Omega^\text{RVE}$. This strategy creates a high-order periodic basis that satisfies standard periodic conditions as stated in Section 3.

At the implementation level, periodicity can be easily treated as a linear constraint on the approximation space during or after the assembly stage, where each basis function $B_\toVect{i}(\xi,\eta,\tau)$ and their corresponding images are added to form the periodic basis function $B_\toVect{i}^\text{P}(\xi,\eta,\tau)$.

\begin{remark}[\it Preventing rigid body translation in periodic simulations]A subtle implementation detail arises in the situation where periodicity is enforced along all the spatial dimensions. Then, the solution of boundary value problems is determined up to rigid body translations, and therefore the resulting system matrix is singular. This issue is easily fixed by setting an arbitrary value of any degree of freedom for each approximated field (\eg enforcing the first degree of freedom of each field to zero). Recall that, in the current framework, rigid body rotations are always prevented, i.e. by enforcing $\overline{\textbf{W}}=0$, as explained in Section \ref{Se:MacroscopicKinematics}.
\end{remark}

\subsection{High-order generalized-periodic approximation space}
\label{Se:GeneralizedPeriodicBasis}

\begin{figure}[htb!]
\centering
    \includegraphics[width=0.6\textwidth]{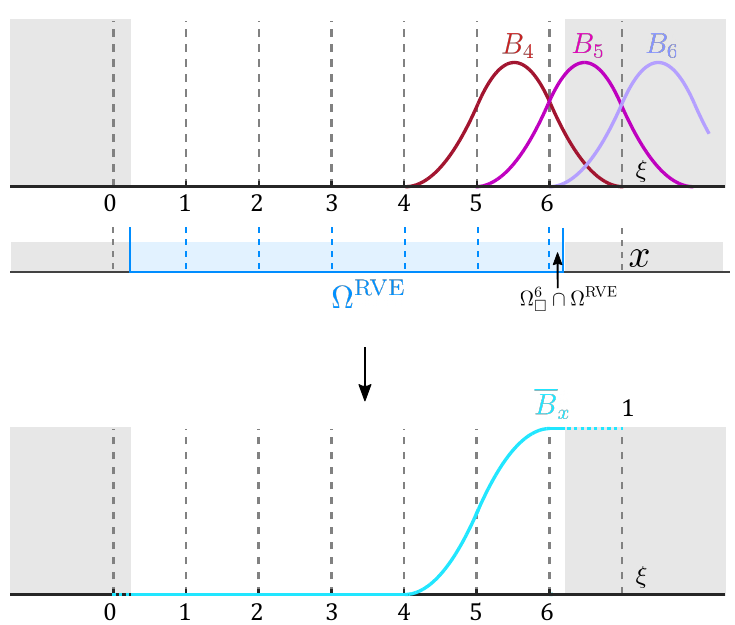}
    \caption[Univariate generalized periodic basis of degree $p=2$]
    {Univariate generalized periodic basis of degree $q=2$. Originally (top), the functional space over the cut cell $\Omega^6_\square\cap\Omega^\text{RVE}$ is spanned by three B-spline basis functions $B_i(\xi),i=\{4,5,6\}$. The addition of these basis functions yields the global basis function $\overline{B}_x(\xi)$ which inherits the regularity of the original B-spline basis and fulfills the admissibility condition in \eq\eqref{eq_generalizedcondition}. The union of $\overline{B}_x$ with the periodic B-spline bases in \fig\ref{fig_PeriodicBasis} spans a high-order generalized-periodic functional space on $\Omega^\text{RVE}$.}
\label{fig_GeneralizedPeriodicBasis}
\end{figure}

By virtue of \eq \eqref{eq:StateDecomposition}, a generalized-periodic approximation space can be constructed by complementing a periodic space with a functional space spanned by global basis functions $\overline{B}_\zeta(\xi,\eta,\tau)$ that fulfill 
\begin{subequations}\label{eq_generalizedcondition}\begin{align}
    \overline{B}_\zeta\left(\rho=n_\rho\right)-\overline{B}_\zeta\left(\rho=0\right)&=\delta_{\zeta \varrho},\\
    \partial^n \overline{B}_\zeta\left(\rho=n_\rho\right)-\partial^n \overline{B}_\zeta\left(\rho=0\right)&=0,
\end{align}\end{subequations}
where $\zeta,\varrho\in \{\xi,\eta,\tau\}$ and $\delta_{\zeta\varrho}$ is the Kronecker delta. 
Since the decomposition in \eq\eqref{eq:StateDecomposition} is not unique, the basis function $\overline{B}_\zeta(\xi,\eta,\tau)$ is neither unique. The simplest definition for an admissible $\overline{B}_\zeta(\xi,\eta,\tau)$ is perhaps the linear function $\overline{B}_\zeta(\xi,\eta,\tau)=\zeta/n_\zeta$. However, it spans throughout the whole $\Omega^\text{RVE}$, substantially increasing the fill-in of the resulting system matrix. An efficient alternative in the context of B-spline bases, that involves minimal fill-in and a very easy implementation, consists on defining $\overline{B}_\zeta(\xi,\eta,\tau)$ as the addition of all the non-vanishing original B-spline bases $B_i(\xi,\eta,\tau)$ on the cell $\Omega_\square^c$ intersected by $\partial\Omega^\text{RVE}=L_\zeta$, as illustrated in \fig\ref{fig_GeneralizedPeriodicBasis}. Thanks to the partition of unity property of B-spline bases, $\overline{B}_\zeta(\xi,\eta,\tau)$ evaluates to 1 within the aforementioned cell and 0 for all derivatives. At the opposite boundary $\partial\Omega^\text{RVE}=0$, $\overline{B}_\zeta(\xi,\eta,\tau)$ and its derivatives vanish, fulfilling the conditions in \eq\eqref{eq_generalizedcondition}. 
The generalized-periodic functional space spanned by $\overline{B}$ and $B^\text{P}$ inherits the regularity of the original B-spline functional space in Section \ref{Se:ImmersedBoundaryApproach}, that is, $C^q-1$ continuity.
In the same way as the periodic bases, $\overline{B}_\zeta(\xi,\eta,\tau)$ can also be implemented as a linear constraint on the original approximation space during or after the assembly stage.

In conclusion, the components of $\Displacement$ and $\phi$ are approximated by the  high-order generalized-periodic functional spaces spanned by the  basis functions $\left\{B_i^\text{P}(\xi,\eta,\tau);\overline{B}_\zeta(\xi,\eta,\tau)\right\}$ and control variables $\left\{\Displacement^\text{P}, \phi^\text{P};\overline{\Strain},\overline{\E}\right\}$ as follows:

\begin{subequations}
\begin{align}
    &u_a (x,y,z) \approx \sum_i \left[B^\text{P}_i \circ\varphi^{-1}(x,y,z)\right] u^\text{P}_a +  \left[\overline{B}_b \circ\varphi^{-1}(x,y,z)\right]\overline{\varepsilon}_{ab} \\
    &\phi (x,y,z) \approx \sum_i \left[B^\text{P}_i \circ\varphi^{-1}(x,y,z)\right] \phi^\text{P} +  \left[\overline{B}_b \circ\varphi^{-1}(x,y,z)\right]\overline{E}_{b}
\end{align}
\end{subequations}

\begin{remark}
 The critical basis functions of the high-order generalized-periodic functional space, that is, those whose support is intersected by $\partial\Omega$ in a very small proportion, can be stabilized by means of the extended B-spline stabilization technique explained in \cite{codony2019immersed} which is a further linear constraint on the approximation space, or by any other means as standard in unfitted methods (see for instance the Ghost penalty method \cite{burman2010ghost} or the artificial stiffness approach \cite{Duster2008,schillinger2015finite} among others). 
\end{remark}

\subsection{Enforcement of macroscopic kinematics}
\label{Se:EnforcementMacroscopicKinematics}
The system of equations that results from discretizing the weak form \eqref{eq:FinalWeakFormGeneralized} on the generalized periodic space from Section \ref{Se:GeneralizedPeriodicBasis} is

\begin{equation}
\label{eq:System}
    \left(
    \begin{aligned}
        \bm{K}_{PP} &&& \bm{K}_{PG}\\
        \bm{K}_{GP} &&& \bm{K}_{GG}
    \end{aligned}
    \right)
    \left(
    \begin{aligned}
        \bm{X}_{P} \\
        \bm{X}_{G} 
    \end{aligned}
    \right)
    =
    \left(
    \begin{aligned}
        &\ 0 \\
        &\bm{f}_{G} 
    \end{aligned}
    \right)
\end{equation}
where the subscript $P$ denotes the periodic basis functions stated in Section \ref{Se:PeriodicBasis} and the subscript $G$ denotes the global basis functions described in Section \ref{Se:GeneralizedPeriodicBasis}. In our case, as stated in Section \ref{Se:energyequivalence}, we have $\bm{X}_G=\{\overline{\Strain},\overline{\E}\}$, $\bm{X}_P=\{\Displacement^{\textnormal{P}},\phi^{\textnormal{P}}\}$ and $\bm{f}_G=|\Omega^\text{RVE}|\{\overline{\bm{\sigma}},\overline{\bm{D}}\}$. 
The sets $\bm{X}_G$ and $\bm{f}_G$ are split in two subsets, one corresponding to macroscopic Dirichlet conditions $\bm{X}_G^D$ and $\bm{f}_G^D$, and another one corresponding to macroscopic Neumann conditions $\bm{X}_G^N$ and $\bm{f}_G^N$. Macroscopic Dirichlet conditions \eq \eqref{eq:MacroscopicDirichlet} are enforced strongly on the system of equations \ref{eq:System} by prescribing the values of $\bm{X}_G^D$, and microscopic Neumann conditions \eq \eqref{eq:MacroscopicNeumann} are enforced by prescribing the values of $\bm{f}^N_G$.
\begin{remark}
So far, the macroscopic conditions have been applied along the directions $(x,y,z)$ of the Cartesian frame. However, macroscopic conditions can be applied along a rotated frame by considering
\begin{align}
    &\overline{\Strain}^\bm{R}=\bm{R}\cdot\overline{\Strain}\cdot\bm{R}^T, &&&
    &\overline{\bm{\sigma}}^\bm{R}=\bm{R}\cdot\overline{\bm{\sigma}}\cdot\bm{R}^T,\nonumber \\
    &\overline{\E}^\bm{R}=\bm{R}\cdot\overline{\E}, &&& 
    &\overline{\bm{D}}^\bm{R}=\bm{R}\cdot\overline{\bm{D}}.
\end{align}
where $\bm{R}$ is a given rotation matrix from the Cartesian frame to the rotated frame. This approach is very convenient, specially in sensitivity analysis w.r.t. loading direction (see Section \ref{2DflexoDevice}), since a continuous response can be obtained by continuously incrementing the rotation angle covering all the parameter space. This task can be performed very efficiently by implementing a for loop during or after the assembly stage, avoiding the re-computation of volume integrals in the approximation space. The new rotated system of equation is

\begin{equation}
    \left(
    \begin{aligned}
        \bm{K}_{PP} &&& \bm{K}_{PG}^{\bm{R}}\\
        \bm{K}_{GP}^{\bm{R}} &&& \bm{K}_{GG}^{\bm{R}}
    \end{aligned}
    \right)
    \left(
    \begin{aligned}
        \bm{X}_{P} \\
        \bm{X}_{G}^\bm{R} 
    \end{aligned}
    \right)
    =
    \left(
    \begin{aligned}
        &\ 0 \\
        &\bm{f}_{G}^\bm{R} 
    \end{aligned}
    \right),
\end{equation}
where $\bm{K}_{PG}^{\bm{R}}=\bm{K}_{PG}\cdot\bm{R}^T$, $\bm{K}_{GP}^{\bm{R}}=\bm{R}\cdot\bm{K}_{GP}$, $\bm{K}_{GG}^{\bm{R}}=\bm{R}\cdot\bm{K}_{GG}\cdot\bm{R}^T$, 
$\bm{X}_G^\bm{R}=\{\overline{\Strain}^{\bm{R}},\overline{\E}^{\bm{R}}\}$ and $\bm{f}_G^\bm{R}=|\Omega^\text{RVE}|\{\overline{\bm{\sigma}}^{\bm{R}},\overline{\bm{D}}^{\bm{R}}\}$. The macroscopic Dirichlet and Neumann conditions can be applied directly to $\bm{X}_G^\bm{R}$ and $\bm{f}_G^\bm{R}$ along the directions of the rotated frame.
\end{remark}

\section{Numerical examples}\label{Se:examples}
We present next several examples aimed at illustrating the capabilities of the proposed framework. The first example in Section \ref{Se:ValidationExample} shows a validation of the proposed approach for flexoelectric RVE, where results on a unit cell under macroscopic Dirichlet conditions are compared against those on a replicated periodic structure considering only equivalent Dirichlet conditions on their actual edges. The second example in Section \ref{2DflexoDevice} is an application of a 2D flexoelectric device, where we compare the performance considering stress-free or strain-free conditions. Also, we simulate different loading directions in order to study its anisotropic behaviour. In the last example in Section \ref{Se:3Dflexoelectric} we show a 3D flexoelectric device and we compare its response against a 2D slice under plain strain conditions. The material tensors needed for this Section are explained in Appendix A.    
\subsection{Validation}
\label{Se:ValidationExample}

For validation purposes, in this Section we analyze the response of a large periodic structure under homogeneous macroscopic vertical deformation directly and compare the results with those obtained on the RVE with generalized periodic conditions. The structure is a periodic concatenation of squares with equilateral triangular voids along the $x$ and $y$ directions. The length of the side of the square is $\SI{4}{\micro\meter}$ and the triangle has a side length of $\frac{3\sqrt{3}}{2}$  $\SI{}{\micro\meter}$. This structure proposed in \cite{Sharma2010} has been shown to mobilize and accumulate the flexoelectric effect to produce a macroscopic electric response under macroscopic homogeneous deformation, \ie it behaves as an apparent piezoelectric, for any dielectric base material \cite{mocci2021geometrically}.

A macroscopic vertical compressive stran is imposed, 
\begin{align}
&\overline{\varepsilon}_{yy}=-0.1,
\label{bc1}
\end{align}
while it is allowed to deform freely in the other directions. 

and all other components of the macroscopic strain and the macroscopic electric field are let free. This means that those components of the macroscopic stress and the macroscopic electric displacement are set to $0$: 

\begin{equation}
    \overline{\sigma}_{xx} =
    \overline{\sigma}_{xy} =
    \overline{\sigma}_{yx} =
    \overline{D}_{x} = \overline{D}_{y} = 0.
\end{equation}
This condition represents imposing an homogeneous unconfined compression along the $y$-direction. The material properties are reported in Table \ref{Tabla1}.

\begin{figure}[tb!]
\centering
    \includegraphics[width=0.9\textwidth]{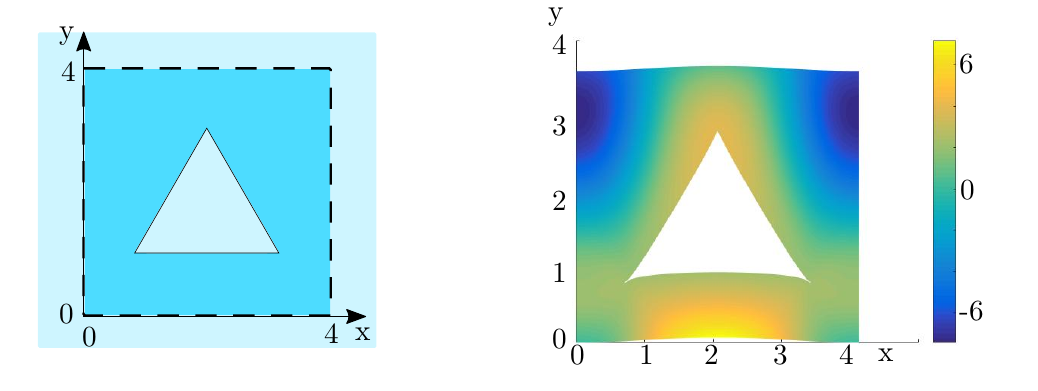}
    \caption{(Left) Unit cell simulated. (Right) Deformed shape and electric potential distribution inside a unit cell, considering generalized periodicity.}
\label{fig_UnitCellEx_5.1}
\end{figure}

\begin{table}[!htp]
	\centering
	\caption{Material in Section \ref{Se:ValidationExample}}
	\begin{tabular}{cccccccc}
		\hline
		\hline
		   E 		& 
		   $\nu$ & 
		   $\ell_{\rm mech}$ &
		   $\kappa$ &
		   $\ell_{\rm elec}$ &
		   $\mu_L$ &
		   $\mu_{T}$ &
		   $\mu_{S}$ \\
		   $[GPa]$ 	&
		   - 		&
		   [$n$m]			  &
		   [$nC/V$m] &
		   [$n$m]		   &
		   [$n C/$m] &
		   [$ nC/$m] &
		   [$n C/$m] \\
		\hline
		   152 	& 0.33 & 1 & 45 & 0 & 40 & 40 & 0\\
		\hline
		\hline
	\end{tabular}\label{Tabla1}
\end{table} 
\fig \ref{fig_UnitCellEx_5.1} shows the unit cell and the resulting electric potential distribution. The macroscopic electric field of the RVE is $\overline{E}_y=-1.6031 V/m$. Now, we consider a stack of $N\in\{1,\dots,20\}$ concatenated cells under prescribed displacements on top and bottom faces matching, in the limit of $N\rightarrow\infty$, the previous generalized periodicity conditions, i.e.

\begin{equation}
\Displacement|_{y=4N}=(0,-0.4N) $\si{\micro\meter}$,\quad\quad\quad
\Displacement|_{y=0}=\bm{0}.
\end{equation}

\begin{figure}[htb!]
\centering
    \includegraphics[width=1\textwidth]{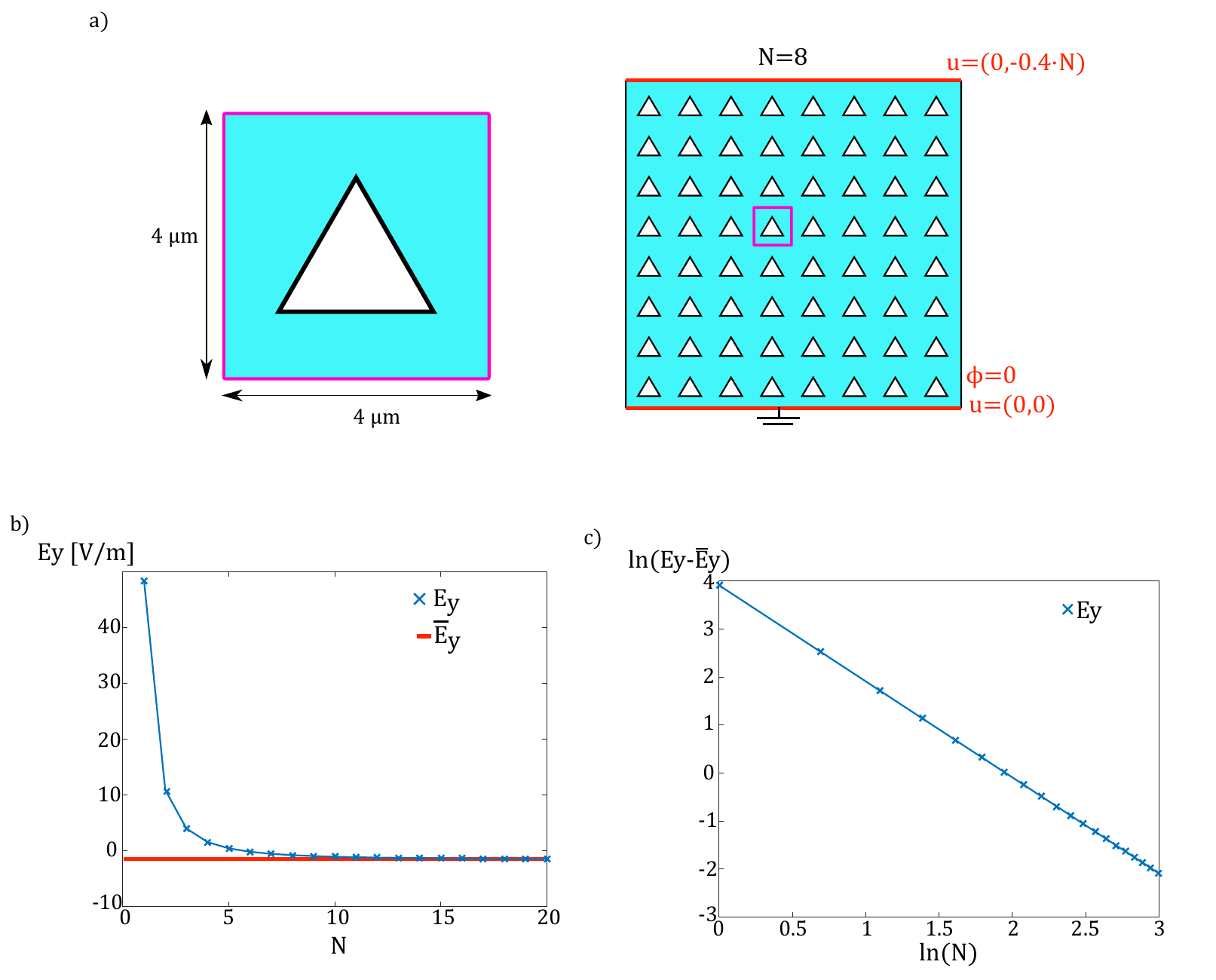}
    \caption{a) Structure formed by $N=8$ cells per dimension. The microscopic Dirichlet boundary conditions are depicted in red and microscopic homogeneous Neumann conditions are applied on all other boundaries. b) Plot of the electric field resulting from stacking $N$ cells per dimension versus the number of cells stacked. The red line is the value resulting considering generalized periodicity conditions. c) Plot of the difference between the electric field resulting from stacking $N$ cells per dimension and the one obtained using generalized periodicity conditions, versus the number of cells stacked. }
\label{fig_NcellEx_5.1}
\end{figure}

The bottom face is electrically grounded ($\phi=0$) and we compute the electric potential drop between the top and bottom boundaries ($\Delta\phi$) as the mean of the electric potential at the top face. Homogeneous Neumann boundary conditions are applied in all other boundaries.  \fig \ref{fig_NcellEx_5.1} shows the periodic structure and the boundary conditions applied for $N=8$.

For stacks of a large enough number of unit cells, we expect that the difference of electric potential between top and bottom faces divided by the vertical length tends to minus the macroscopic electric field $\overline{E}_y$, i.e. 
\begin{equation}
\lim\limits_{N\rightarrow\infty}\frac{-\Delta\phi}{L_y}=\lim\limits_{N\rightarrow\infty}E_y=
\overline{E}_y.
\label{limit}
\end{equation}

For quantitative validation purposes, we plot in \fig\ref{fig_NcellEx_5.1} the value of $E_y$ for stacks of increasing number of unit cells, from $N=1$ to $N=20$ and compare it against the macroscopic electric field $\overline{E}_y$ obtained in the generalized periodic unit cell simulation, \fig\ref{fig_UnitCellEx_5.1}.  From this plot, the limit in \eq\eqref{limit} is apparent. The electromechanical behaviour of a unit cell under generalized periodicity conditions is representative of the behaviour of a unit cell in the bulk of a periodic structure. The difference between the value considering $N=20$ cells and the one obtained using generalized periodicity is less than $0.05\%$. 
This validation is extensible to any other set of generalized periodic conditions, resulting in the same conclusion.

\subsection{2D flexoelectric metamaterial}\label{2DflexoDevice}
\begin{figure}[htb!]
\centering
    \includegraphics[width=1\textwidth]{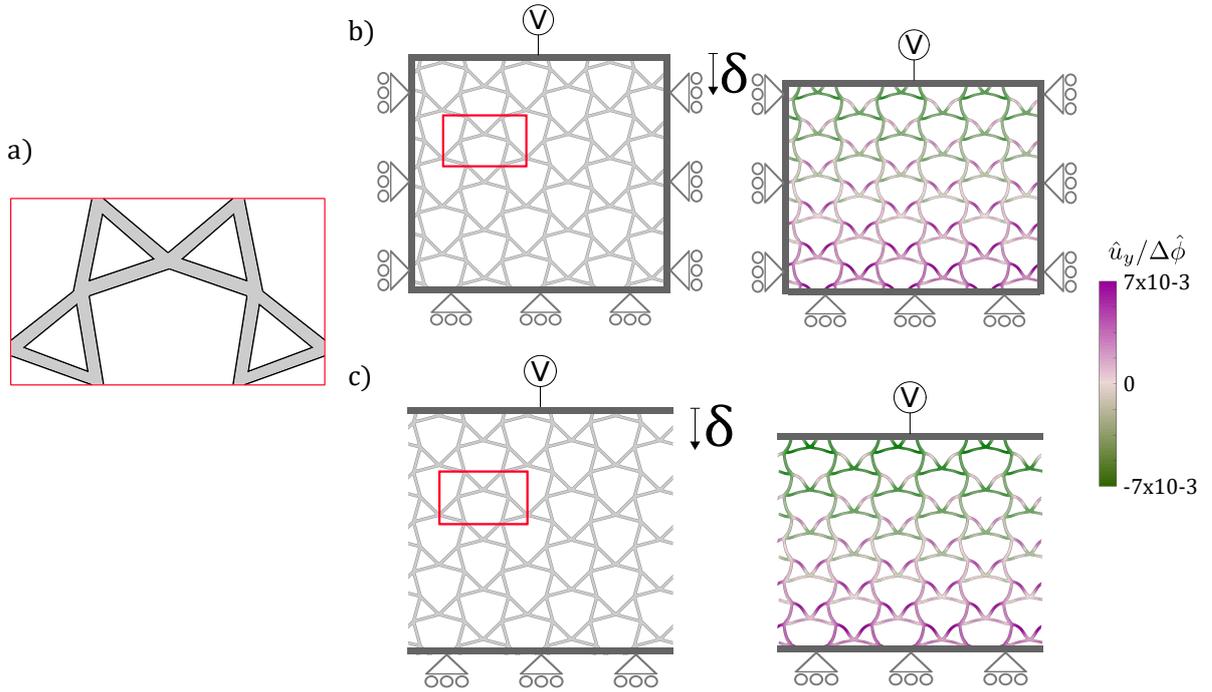}
    \caption{Geometrically-polarized, bending-dominated lattice in actuation mode. The RVA is represented in a). The lattice is resolved with respect to the vertical displacements upon a macroscopic electric field in the vertical direction. Horizontally, we considered b) classic periodicity \ie a strain-free condition $\overline{\strain}_{xx}=0$ and c) stress-free condition $\overline{\sigma}_{xx}=0$. The normalized y-displacements $\hat{u}_y=u_y/\ell$ are depicted on the deformed configuration considering the normalized electric potential $\hat{\phi}=\phi\kappa/\mu$.}
\label{fig_NumericalEx_5.2}
\end{figure}

Flexoelectricity might represent an alternative route for technologies based on electromechanical transduction. If proper designed, architected metamaterials can endow any dielectric with apparent piezoelectricity \cite{sharma2007possibility,mocci2021geometrically}. 
Here, we consider a low-area fraction, bending-dominated  metamaterial, as proposed by \cite{mocci2021geometrically}. The lattice behaves as an apparent piezoelectric due to the flexoelectric effect mobilized in the micro-constituents and their non-centrosymmetric arrangement which avoids internal cancellations and accumulates the response. The lattice is geometrically constructed by filling the two-dimensional space with stacked periodic unit cells characterized by thin features having thickness $t=160$ nm and length $\ell=1.6\upmu$m.  To efficiently evaluate the response of the flexoelectric lattice we reduce the computational domain to a RVE, highlighted in red in \fig\ref{fig_NumericalEx_5.2}(a) , under generalized periodicity conditions. Figs.~\ref{fig_NumericalEx_5.2}(b,c) show two different setups for an actuator mode. In both configurations a homogeneous vertical electric field $\overline{E}_{y}$ is applied whereas mechanically the lattice is free to deform vertically and $\overline{\varepsilon} _{yy}$ is computed. Horizontally, we consider two limit cases: (1) a \textit{rigid} device by considering a strain-free configuration thus imposing classic periodicity $\overline{\varepsilon}_{xx}=\overline{\varepsilon}_{xy}=0$  (cf.~\fig\ref{fig_NumericalEx_5.2}b), and (2) a \textit{soft} configuration, \ie stress-free condition $\overline{\sigma}_{xx}=\overline{\sigma}_{xy}=0$ (cf.~\fig\ref{fig_NumericalEx_5.2}c). The material parameters of the base material are chosen accordingly to \cite{mocci2021geometrically} and reported in table \ref{table:Material Parameters}. 

Although qualitatively the responses do not differ for the two setups, two different apparent piezoelectric coefficients $\hat{d}$, $\bar{d}$ are quantified, as
\begin{align}
	\quad \hat{d} = \frac{\partial \bar{{\varepsilon}}}{\partial \bar{{E}}}\biggl\lvert_{\bar{\epsilon}_{xx}=0}=0.51 \ pm/V,
	\quad \bar{d} = \frac{\partial \bar{{\varepsilon}}}{\partial \bar{{E}}}\biggl\lvert_{\bar{\sigma}_{xx}=0}=0.42 \ pm/V, 
 \label{eq_coefs}
\end{align}
where the rigid device shows over a $20\%$ performance improvement with respect to the soft device. 
\begin{table}[!tb]
	\centering
	\caption{Material in Section \ref{2DflexoDevice}}
	\begin{tabular}{c|cccccccc}
		\hline
		\hline
		Material & Y 		& $\nu$ & $\ell_{\rm mech}$ & $\kappa$ & $\ell_{\rm elec}$ &$\mu_{L}$ & $\mu_{T}$ & $\mu_{S}$ \\
		&[$GPa$] 	&- 		&[$n$m]			  &[$nC/V$m] &[$n$m]		   &[$\mu C/$m] &[$\mu C/$m] &[$\mu C/$m] \\
		\hline
		BST & 152 	& 0.33 & 50 & 8 & 300 & 1.21 & 1.10 & 0.055\\
		\hline
		\hline
	\end{tabular}\label{table:Material Parameters}
\end{table} 

\begin{figure}[htb!]
\centering
    \includegraphics[width=1\textwidth]{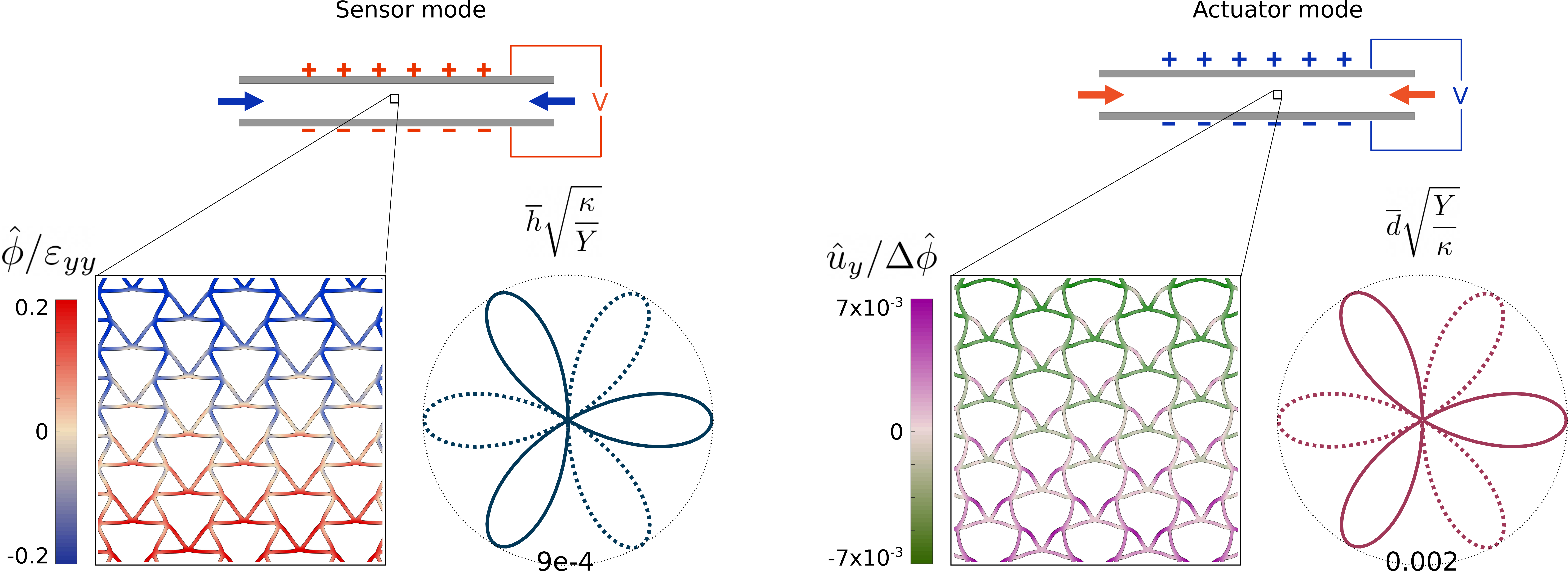}
    \caption{Anisotropy of the normalized apparent piezoelectric coefficients $\bar{h}\sqrt{\kappa/Y}$ and $\bar{d}\sqrt{Y/\kappa}$ for sensor and actuator mode, respectively. In the polar plots, solid lines are used to depict positive values of the apparent piezoelectric coefficients, while dashed lines are used to indicate negative values. The normalized electric potential $\hat{\phi}=\phi\kappa/\mu$ and y-displacements $\hat{u}_y=u_y/l$ are plotted on the deformed configuration.}
\label{fig_NumericalEx_Anisotropy}
\end{figure}

The anisotropy of the architected lattice can also be studied by considering different orientations of
the mechanical and electrical loadings (see remark in Section \ref{Se:EnforcementMacroscopicKinematics}). The response of the lattice is represented by the apparent piezoelectric coefficients $\bar{h}$ and $\bar{d}$ suitably normalized with respect to the nominal Young's modulus $Y$ and dielectric permittivity $\kappa$ of the base material. Two polar plots are obtained under the applied rotated strain and electric field, respectively for sensor and actuator mode, cf.~\fig\ref{fig_NumericalEx_Anisotropy}. As expected the behavior is highly anisotropic. No polarization is observed when loaded horizontally as the lattice is not geometrically polarized along the horizontal direction. The lattice exhibits symmetry with respect to $2\pi/3$ rotations and one planar mirror symmetry, accurately recovered in the polar plots, where solid and dotted lines represent positive and negative values of the normalized apparent piezoelectric coefficient, respectively. 

\subsection{3D flexoelectric metamaterial}
\label{Se:3Dflexoelectric}
\begin{figure}[htb!]
\centering
    \includegraphics[width=1\textwidth]{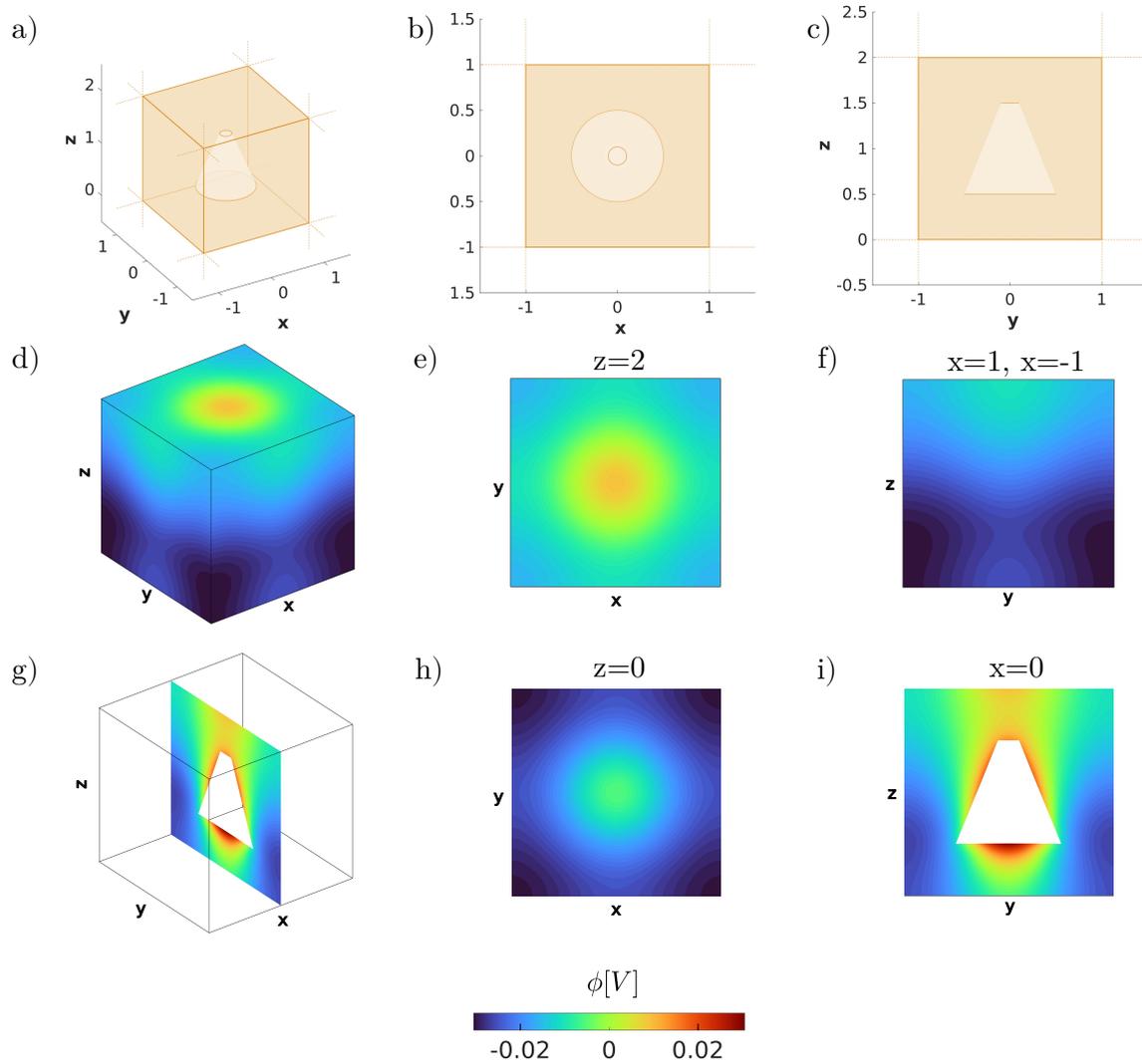}
    \caption{a-c) Geometrical model and d-i) electric potential distribution in different cross sections, upon macroscopic deformation in a 3D RVE with geometrically-polarized truncated conical void embedded in a flexoelectric matrix.}
\label{fig_NumericalEx_5.4}
\end{figure}

The concept of an architected dielectric endowed with apparent piezoelectricity can naturally be extended to 3D. Here we consider an architected material consisting of voids shaped as truncated cones of radius $1\upmu$m and $0.2\upmu$m and height $1\upmu$m embedded in a flexoelectric matrix (material properties in Table \ref{table:Material Parameters}). The periodic RVE (cf.~Figs.~\ref{fig_NumericalEx_5.4}a-c), where generalized periodic boundary conditions are enforced, is a cube having size $\ell=2\upmu$m. Figs.~\ref{fig_NumericalEx_5.4}(d-i) show the electric potential distribution $\phi$ in the RVE in sensing mode, upon macroscopic deformation $\overline{\varepsilon}_{zz}=-0.1$, whereas standard periodicity is applied for the solution fields in the other principal directions $\overline{\varepsilon}_{xx}=\overline{\varepsilon}_{yy}=\overline{\varepsilon}_{xz}=\overline{\varepsilon}_{xy}=\overline{\varepsilon}_{yz}=0$ and $\overline{D}_x=\overline{D}_y=\overline{D}_z=0$. Similarly to the electromechanical response showed in Section \ref{2DflexoDevice}, the geometrically-polarized inclusion produces a macroscopic electric field $\overline{E}_{z}$, as clearly showed in the different sections of the unit cell, reported in Figs.~\ref{fig_NumericalEx_5.4}(d-i). 

\begin{figure}[htb!]
\centering
    \includegraphics[width=0.7\textwidth]{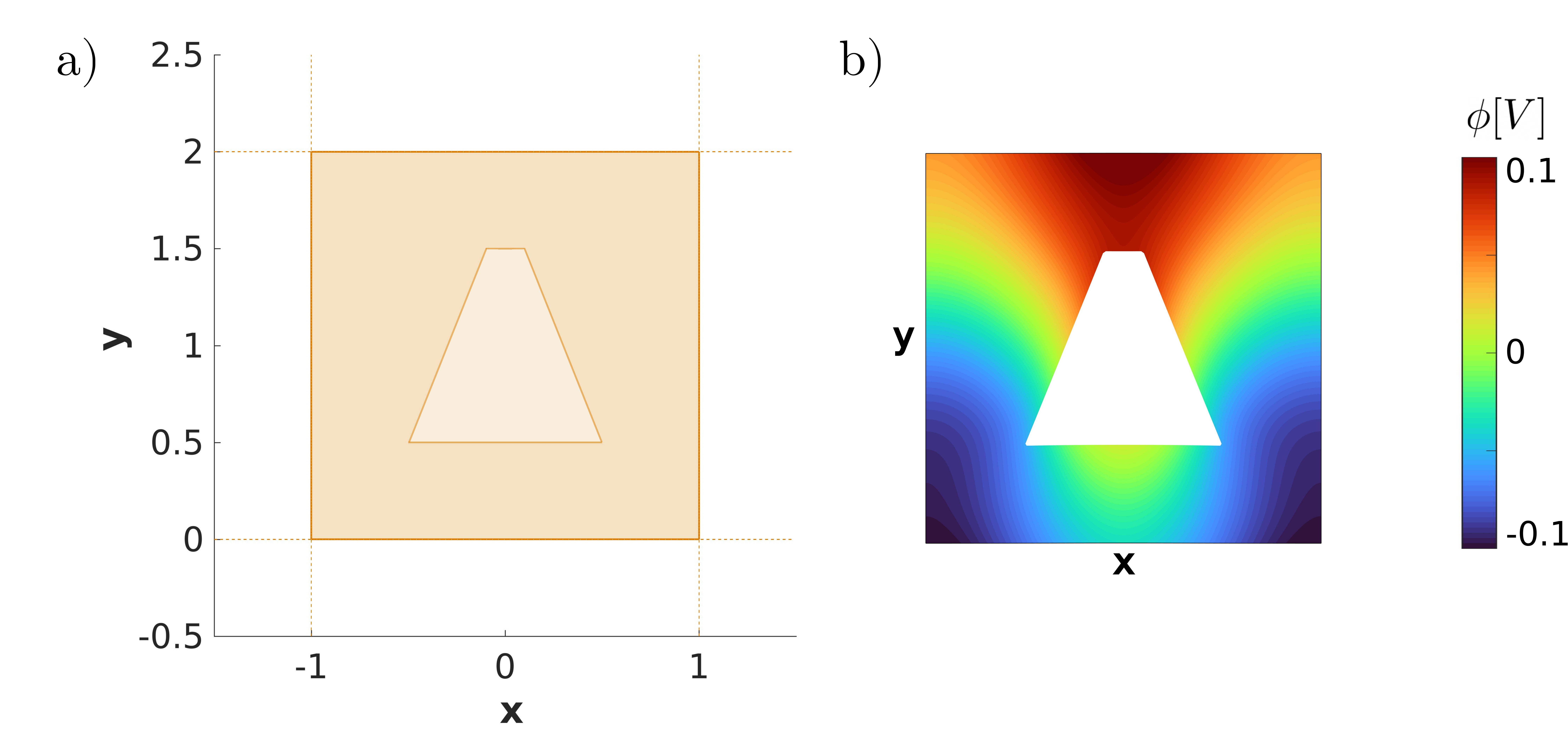}
    \caption{a) Geometrical model and b) electric potential distribution upon homogeneous deformation in a 2D unit cell with geometrically-polarized truncated conical hole embedded in a flexoelectric matrix.}
\label{2D_SquarePyramid}
\end{figure}

It is also instructive to compare the 3D electromechanical response with a similar 2D RVE considering plane strain. In this regard, we reproduce a 2D unit cell (cf.~\fig\ref{2D_SquarePyramid}) as a replica of the cross section $x=0$ of the 3D structure in \fig\ref{fig_NumericalEx_5.4}. Similarly, we enforce a macroscopic deformation $\overline{\varepsilon}_{yy}=-0.1$ with $\overline{\varepsilon}_{xy}=\overline{\varepsilon}_{xx}=0$ and $\overline{D}_{x}=\overline{D}_{y}=0$ and we resolve the macroscopic electric field along the vertical direction $\overline{E}_{y}$, \fig\ref{2D_SquarePyramid}(b). The difference in the solution field in \fig\ref{fig_NumericalEx_5.4}(i) and \fig\ref{2D_SquarePyramid}(b) highlights that 3D simulations are required to comprehensively quantify the response. 

\section{Conclusions}\label{Se:Conclusions}

We have proposed a simple and elegant method to computationally account for the high-order generalized periodicity conditions that arise when analyzing a periodic metamaterial with a RVE. Such computational homogenization problem arises when the base material is modeled mathematically with higher-order partial differential equations. Here, we have focused on architected materials for electromechanical transduction exhibiting apparent piezoelectricity and made from a non-piezoelectric and flexoelectric base-material. The generalized periodicity conditions enable us to impose macroscopic fields (strain/stress, electric field/electric displacement) and recover the dual response, and hence to characterize the tensorial mechanical, electrical, and piezoelectric response.

The proposed method constructs a high-order generalized-periodic approximation space for the state variables, as an extension of a B-spline approximation over unfitted Cartesian grids. Dirichlet and Neumann macroscopic conditions are enforced strongly along arbitrary directions. This is particularly useful to study material anisotropy without the need tor rotating the RVE. The method is verified by direct comparison against an increasingly large architected material. Several illustrative examples in 2D in 3D show the potential of the method, including the systematic study of anisotropic apparent piezoelectricity  of non-piezoelectric lattices. Although we focus here on flexoelectricity, the formulation is readily extensible to other high-order boundary value problems. Being devoid of the stability issues or penalty parameters of methods imposing the generalized periodicity constraints in the variational formulation, our method is computationally robust and can be integrated in topology optimization frameworks at the RVE level \cite{Greco2023}.

\section*{Acknowledgments}
This work was supported by the Generalitat de Catalunya  (``ICREA Academia'' award for excellence in research to I.A., and Grant No.~2017-SGR-1278), the European Research Council (StG-679451 to I.A.), the FI-AGAUR grant and the grant CEX2018-000797-S funded by MCIN/AEI/10.13039/501100011033 (Severo Ochoa Centre of Excellence 2019-2023 to CIMNE).  D.C.~acknowledges the support of the Spanish Ministry of Universities through the Margarita Salas fellowship (European Union-NextGenerationEU).

\appendix

\section{Material Tensors}

In the following appendix, the material tensors are defined, as in \cite{codony2019immersed,codony2021mathematical}. They are described component-wise (non-zero components) and $d$ is the number of dimensions of the physical space.

We use an isotropic elasticity tensor defined in terms of the Young modulus $E$ and Poisson ratio $\nu$ as 
\begin{align}
\elast_{iiii}&=C_L,&&& i&=1,\dots,d;\nonumber\\
\elast_{iijj}&=C_T,&&& i,j&=1,\dots,d\ : \ i\neq j;\nonumber\\
\elast_{ijij}=\elast_{ijji}&=C_S,&&& i,j&=1,\dots,d\ : \ i\neq j,
\end{align}
where the parameters $C_L$, $C_S$ and $C_T$ are 
\begin{align}\label{eq_c3}
C_L=\frac{E\left(1-\nu\right)}{(1+\nu)(1-2\nu)},&&
C_T=\frac{E\nu}{(1+\nu)(1-2\nu)},&&
C_S=\frac{C_L-C_T}{2}=\frac{E}{2(1+\nu)}.
\end{align}

We use a sixth-order tensor to describe strain gradient elasticity. We consider an isotropic version of the general model in \cite{Mindlin1968a} that is described in  \cite{Altan1997}. The strain gradient tensor depends on the Young modulus $E$, the Poisson ratio $\nu$ and the mechanical length scale $\ell_{\text{mech}}$ as
\begin{align}
\strGr_{iikiik}&=\ell_{\text{mech}}^2C_L,&&& i,k&=1,\dots,d;\nonumber\\
\strGr_{iikjjk}&=\ell_{\text{mech}}^2C_T,&&& i,j,k&=1,\dots,d\ : \  i\neq j;\nonumber\\
\strGr_{ijkijk}=\strGr_{ijkjik}&=\ell_{\text{mech}}^2C_S,&&& i,j,k&=1,\dots,d\ : \  i\neq j
\end{align}
where the parameters $C_L$, $C_S$ and $C_T$ are defined in \eq\eqref{eq_c3}.

We use a second-order tensor to describe isotropic dielectricity, which depends on a parameter $\epsilon$ as 
\begin{align}
\epsilon_{ii}&=\epsilon,&&& i&=1,\dots,d.
\end{align}

Isotropic gradient dielectricity is represented by the fourth-order tensor $M$. We take a simple form depending on the electric permittivity $\epsilon$ and the dielectric length scale $\ell_\text{elec}$ as 
\begin{align}
	&M_{ijij}=\epsilon\ell_\text{elec}^2 &&& &i,j=1,\dots,d.
\end{align}
Piezoelectricity is represented by the third-order tensor $\Piezo$, where tetragonal symmetry is considered,
which has a principal direction. It involves longitudinal, transversal and shear couplings represented by the parameters $\piezo_L$, $\piezo_T$ and $\piezo_S$, respectively. For a material with principal direction $\x_1$, the piezoelectric tensor $\Piezo^{<\x_1>}$ is
\begin{align}
{\piezo^{<\x_1>}}_{111}&=\piezo_L;\nonumber\\
{\piezo^{<\x_1>}}_{1jj}&=\piezo_T,&&& j&=2,\dots,d;\nonumber\\
{\piezo^{<\x_1>}}_{j1j}={\piezo^{<\x_1>}}_{jj1}&=\piezo_S,&&& j&=2,\dots,d.
\end{align}
The piezoelectric tensor $\Piezo$ oriented in an arbitrary direction $\toVect{d}$ is obtained by rotating $\Piezo^{<\x_1>}$.

Flexoelectricity is represented by the fourth-order tensor $\Flexo$ where cubic symmetry is considered. It leads to a tensor involving longitudinal, transversal and shear couplings represented by the parameters $\flexo_L$, $\flexo_T$ and $\flexo_S$, respectively. The components of the flexoelectric tensor $\Flexo^{<{\x}>}$ of a material oriented in the Cartesian axes are the following:
\begin{align}
{\flexo^{<{\x}>}}_{iiii}&=\flexo_L,&&& i&=1,\dots,d\nonumber;\\
{\flexo^{<{\x}>}}_{ijji}&=\flexo_T,&&& i,j&=1,\dots,d\ : \  i\neq j;\nonumber\\
{\flexo^{<{\x}>}}_{iijj}={\flexo^{<{\x}>}}_{ijij}&=\flexo_S,&&& i,j&=1,\dots,d\ : \  i\neq j.
\end{align}
The flexoelectric tensor $\Flexo$ oriented in an arbitrary orthonormal basis is obtained by rotating $\Flexo^{<\x>}$.

\section{Derivation of alternative expression for macroscopic stress}
This Appendix proofs that the expression of the macroscopic stress in \eq \eqref{eq:DefMacrostress2} is equivalent to that of \eq (18) in Balcells et al.~\cite{Balcells}, that is,

\begin{subequations}
    \begin{align}
        \overline{\sigma}_{i1}=\frac{F_i^x}{L_y}=\frac{1}{L_y}\left(\int_{\Gamma^{L_x}} t_i^{L_x}\dd \Gamma + \sum_{C\in C^{L_x}  }j_i^{L_x} \right), \\
        \overline{\sigma}_{i2}=\frac{F_i^y}{L_x}=\frac{1}{L_x}\left( \int_{\Gamma^{L_y}} t_i^{L_y} \dd \Gamma+\sum_{C\in C^{L_y}  }j_i^{L_y} \right).
    \end{align}
\end{subequations}
We start from the expression derived in \eqref{eq:DefMacrostress2}:
\begin{equation}
    |\Omega^\text{RVE}|\overline{\sigma}_{ij}=\int_\Omega \cauchyStress_{ij}\dd\Omega=\int_\Omega\left( \frac{\partial X_j}{\partial X_k}\cauchyStress_{ik}\right)\dd\Omega=-\int_\Omega\left(X_j\frac{\partial \cauchyStress_{ik}}{\partial X_k}\right)\dd\Omega + \int_\Gamma X_j\cauchyStress_{ik}n_k\dd\Gamma,
\end{equation}
where we have used integration by parts and the divergence theorem. Now we consider the strong form of the problem \eqref{eq_EulerLagrange} with $b_i=0$ and we get
\begin{equation}
    -\int_\Omega\left(X_j\frac{\partial \cauchyStress_{ik}}{\partial X_k}\right)\dd\Omega + \int_\Gamma X_j\cauchyStress_{ik}n_k\dd\Gamma=-\int_\Omega\left(X_j\frac{\partial^2 \hyperStress_{ikq}}{\partial X_k\partial X_q}\right)\dd\Omega + \int_\Gamma X_j\cauchyStress_{ik}n_k\dd\Gamma.
\end{equation}
Considering again integration by parts and the divergence theorem, as done in \cite{codony2021mathematical} to find the strong form of the boundary value problem ,we get

\begin{align}
    &-\int_\Omega\left(X_j\frac{\partial^2 \hyperStress_{ikq}}{\partial X_k\partial X_q}\right)\dd\Omega + \int_\Gamma X_j\cauchyStress_{ik}n_k\dd\Gamma \nonumber =
    -\int_{\Omega}\frac{\partial}{\partial X_k}\left( X_j \frac{\partial \hyperStress_{ikq}}{\partial X_q}\right) \dd\Omega + \int_{\Omega}\frac{\partial X_j}{\partial X_k}\frac{\partial \hyperStress_{ikq}}{\partial X_q}\dd \Omega + \int_\Gamma X_j\cauchyStress_{ik}n_k\dd\Gamma=\\
    &-\int_{\Gamma} X_j \frac{\partial \hyperStress_{ikq}}{\partial X_q}n_k \dd\Gamma + \int_{\Omega}\delta_{jk}\frac{\partial \hyperStress_{ikq}}{\partial X_q}\dd \Omega + \int_\Gamma X_j\cauchyStress_{ik}n_k\dd\Gamma =
    \int_{\Gamma} X_j \left( \cauchyStress_{ik}- \frac{\partial \hyperStress_{ikq}}{\partial X_q}\right) n_k \dd\Gamma + \int_{\Omega}\frac{\partial \hyperStress_{ijq}}{\partial X_q}\dd \Omega= \nonumber \\
    &\int_{\Gamma} X_j \left( \cauchyStress_{ik}- \frac{\partial \hyperStress_{ikq}}{\partial X_q}\right) n_k \dd\Gamma + \int_{\Gamma} \hyperStress_{ijq}n_q\dd \Gamma=
    \int_{\Gamma} X_j \left( \cauchyStress_{ik}- \frac{\partial \hyperStress_{ikq}}{\partial X_q}\right) n_k \dd\Gamma + \int_{\Gamma}\frac{\partial X_j}{\partial X_k}\hyperStress_{ikl} n_l \dd \Gamma.
\end{align}
Now splitting the term $\partial X_j/\partial X_k$ in the normal part and tangential part and applying the surface divergence theorem we arrive to the last expression:
\begin{align}
    &\int_{\Gamma} X_j \left( \cauchyStress_{ik}- \frac{\partial \hyperStress_{ikq}}{\partial X_q}\right) n_k \dd\Gamma + \int_{\Gamma}\frac{\partial X_j}{\partial X_k}\hyperStress_{ikl} n_l \dd \Gamma=\int_{\Gamma} X_j \left( \cauchyStress_{ik}- \frac{\partial \hyperStress_{ikq}}{\partial X_q}\right) n_k \dd\Gamma +
    \int_{\Gamma}\hyperStress_{ikl}n_l\left(\nabla_k^S+n_k\partial ^n  \right) X_j  \dd \Gamma= \nonumber \\
    &\int_{\Gamma} X_j \left( \cauchyStress_{ik}- \frac{\partial \hyperStress_{ikq}}{\partial X_q}\right) n_k \dd\Gamma +
    \int_\Gamma \partial ^n\left(X_j\right)\hyperStress_{ikl}n_ln_k\dd \Gamma+\int_\Gamma\nabla_k^S\left(X_j\right) \hyperStress_{ikl}n_l\dd\Gamma\nonumber =\\
    &\int_{\Gamma} X_j \left( \cauchyStress_{ik}- \frac{\partial \hyperStress_{ikq}}{\partial X_q}\right) n_k \dd\Gamma +
    \int_\Gamma \partial ^n\left(X_j\right)\hyperStress_{ikl}n_ln_k\dd \Gamma-\int_\Gamma\nabla_k^S\left( \hyperStress_{ikl}n_l\right)X_j\dd\Gamma+\int_\Gamma\nabla_k^S\left(\hyperStress_{ikl}n_lX_j\right)\dd\Gamma\nonumber =\\
    &\int_{\Gamma} X_j \left( \cauchyStress_{ik}- \frac{\partial \hyperStress_{ikq}}{\partial X_q}\right) n_k \dd\Gamma +
    \int_\Gamma \partial ^n\left(X_j\right)\hyperStress_{ikl}n_ln_k\dd \Gamma-\int_\Gamma\nabla_k^S\left( \hyperStress_{ikl}n_l\right)X_j\dd\Gamma+\nonumber \\ &\int_\Gamma\nabla_q^S\left(n_q\right)\hyperStress_{ikl}n_kn_lX_j\dd\Gamma+\int_C\jump{\hyperStress_{ikl}n_lm_k}X_j\dd S=\nonumber\\
    &\int_{\Gamma} X_j \left( \left(\cauchyStress_{ik}- \frac{\partial \hyperStress_{ikq}}{\partial X_q}+\nabla_q^S\left(n_q\right)\hyperStress_{ikl}n_l\right) n_k-\nabla_k^S\left( \hyperStress_{ikl}n_l\right)\right) \dd\Gamma +
    \int_\Gamma \partial ^n\left(X_j\right)\hyperStress_{ikl}n_ln_k\dd \Gamma+\int_C\jump{\hyperStress_{ikl}n_lm_k}X_j\dd S.
\end{align}
Substituting the definitions of the traction, double traction and force per unit length, we obtain the following expression:
\begin{align}\label{eq:intMacroTraction}
    &\int_{\Gamma} X_j \left( \left(\cauchyStress_{ik}- \frac{\partial \hyperStress_{ikq}}{\partial X_q}+\nabla_q^S\left(n_q\right)\hyperStress_{ikl}n_l\right) n_k-\nabla_k^S\left( \hyperStress_{ikl}n_l\right)\right) \dd\Gamma +
    \int_\Gamma \partial ^n\left(X_j\right)\hyperStress_{ikl}n_ln_k\dd \Gamma+\int_C\jump{\hyperStress_{ikl}n_lm_k}X_j\dd S=\nonumber \\
    &\int_{\Gamma} X_j t_i \dd\Gamma +
    \int_\Gamma \partial ^n\left(X_j\right)r_i\dd \Gamma+\int_C j_iX_j\dd S.
\end{align}

Following the notation of \cite{Balcells},we split the boundary $\Gamma$ in two sets $\Gamma = \Gamma^{\text{actual}} + \Gamma^{\text{fict}}$, where $\Gamma^{\text{actual}}$ corresponds to the actual boundaries in $\Omega$ and $\Gamma^{\text{fict}} = \partial \Omega \cap \partial \Omega^\text{RVE}$ are the fictitious boundaries that arise due to the restriction of $\Omega$ to a single RVE $\Omega^\text{RVE}$. Assuming that homogeneous Neumann is considered in $\Gamma^{\text{actual}}$, just the integral over $\Gamma^{\text{fict}}$ remains in \eqref{eq:intMacroTraction}. Analogously, the integral over $C$ is restricted to an integral over $C^{\text{fict}} = C \cap \Omega^{\text{RVE}}$. Now, we split $\Gamma^{\text{fict}}$ in 3 different sets, $\Gamma^{\text{fict}}=\Gamma^X\cup \Gamma^Y\cup \Gamma^Z$ corresponding to the planes orthogonal to the X, Y and Z directions respectively, and similarly, $C^{\text{fict}}=C^X\cup C^Y\cup C^Z$. We also split each $\Gamma^\zeta$ in two sets, $\Gamma^{0_\zeta}$ and $\Gamma^{L_\zeta}$ with $\zeta\in\{X,Y,Z\}$. We get
\begin{align}
    &\int_{\Gamma} X_j t_i \dd\Gamma +
    \int_\Gamma \partial ^n\left(X_j\right)r_i\dd \Gamma+\int_C j_iX_j\dd S= \int_{\Gamma^{\text{fict}}} X_j t_i \dd\Gamma +
    \int_{\Gamma^{\text{fict}}} \partial ^n\left(X_j\right)r_i\dd \Gamma+\int_{C^{\text{fict}}} j_iX_j\dd S \nonumber = \\
    &\sum_\zeta \left(\int_{\Gamma^{\zeta}} X_j t_i \dd\Gamma +
    \int_{\Gamma^{\zeta}} \partial ^n\left(X_j\right)r_i\dd \Gamma+\int_{C^{\zeta}} j_iX_j\dd S\right) \nonumber = \\
    &\sum_\zeta \left(\int_{\Gamma^{0_\zeta}} X_j t_i \dd\Gamma +
    \int_{\Gamma^{0_\zeta}} \partial ^n\left(X_j\right)r_i\dd \Gamma+\int_{C^{0_\zeta}} j_iX_j\dd S+\int_{\Gamma^{L_\zeta}} X_j t_i \dd\Gamma +
    \int_{\Gamma^{L_\zeta}} \partial ^n\left(X_j\right)r_i\dd \Gamma+\int_{C^{L_\zeta}} j_iX_j\dd S\right).
\end{align}
Using the equilibrium conditions in \eqref{eq:macrotraction}, we obtain
\begin{align}
    &\sum_\zeta \Bigg(\int_{\Gamma^{0_\zeta}} X_j t_i \dd\Gamma +
    \int_{\Gamma^{0_\zeta}} \partial ^n\left(X_j\right)r_i\dd \Gamma+\int_{C^{0_\zeta}} j_iX_j\dd S+\int_{\Gamma^{L_\zeta}} X_j t_i \dd\Gamma +
    \int_{\Gamma^{L_\zeta}} \partial ^n\left(X_j\right)r_i\dd \Gamma+\int_{C^{L_\zeta}} j_iX_j\dd S\Bigg)\nonumber =
    \\
    &\sum_\zeta \left(\int_{\Gamma^{L_\zeta}} t_i^{L_\zeta}\left(X_j^{L_\zeta}-X_j^{0_\zeta}\right) \dd\Gamma +
    \int_{\Gamma^{L_\zeta}}
    r_i^{L_\zeta}\left(\partial^n\left(X_j^{L_\zeta}\right)+\partial^n\left(X_j^{0_\zeta}\right)\right)\dd\Gamma
    \int_{C^{L_\zeta}} j_i^{L_\zeta}\left(X_j^{L_\zeta}-X_j^{0_\zeta}\right)\dd S\right)=\nonumber
    \\
    &\sum_\zeta \left(L_\zeta\int_{\Gamma^{L_\zeta}} t_i^{L_\zeta}\hat{e}^\zeta_j \dd\Gamma +
    L_\zeta\int_{C^{L_\zeta}} j_i^{L_\zeta}\hat{e}^\zeta_j\dd S\right),
\end{align}
where we have used $\partial^n\left(X_j^{L_\zeta}\right)+\partial^n\left(X_j^{0_\zeta}\right)=0$. Therefore

\begin{equation}
    \overline{\sigma}_{ij}=\frac{1}{|\Omega^\text{RVE}|}\sum_\zeta \left(L_\zeta\int_{\Gamma^{L_\zeta}} t_i^{L_\zeta}\hat{e}^\zeta_j \dd\Gamma +
    L_\zeta\int_{C^{L_\zeta}} j_i^{L_\zeta}\hat{e}^\zeta_j\dd S\right),
\end{equation}

which simplifies in 2D as

\begin{subequations}
    \begin{align}
        \overline{\sigma}_{i1}=\frac{1}{L_y}\left(\int_{\Gamma^{L_x}} t_i^{L_x}  \dd \Gamma+ \sum_{C\in C^{L_x}  }j_i^{L_x}\right) =\frac{F_i^x}{L_y}, \\
        \overline{\sigma}_{i2}=\frac{1}{L_x}\left( \int_{\Gamma^{L_y}} t_i^{L_y} \dd \Gamma+\sum_{C\in C^{L_y}  }j_i^{L_y} \right) =\frac{F_i^y}{L_x}.
    \end{align}
\end{subequations}

The expression of the macro electric displacement is derived analogously. This shows that the macroscopic stress, which is the macroscopic volume average of the microscopic Cauchy stress, corresponds as well to the macroscopic surface averages of the microscopic forces per unit area and length in the RVE boundary.


\begin{thebibliography}{10}
\expandafter\ifx\csname url\endcsname\relax
  \def\url#1{\texttt{#1}}\fi
\expandafter\ifx\csname urlprefix\endcsname\relax\def\urlprefix{URL }\fi
\expandafter\ifx\csname href\endcsname\relax
  \def\href#1#2{#2} \def\path#1{#1}\fi

\bibitem{Engetha2006}
N.~Engheta, R.~W. Ziolkowski, Metamaterials: Physics and Engineering
  Explorations, John Wiley \& Sons, Ltd., 2006.

\bibitem{10.1063/1.3490504}
D.~D. Paul, \href{https://doi.org/10.1063/1.3490504}{{ Optical Metamaterials:
  Fundamentals and Applications }}, Physics Today 63~(9) (2010) 57--58.

\bibitem{Bertoldi:2017aa}
K.~Bertoldi, V.~Vitelli, J.~Christensen, M.~van Hecke,
  \href{https://doi.org/10.1038/natrevmats.2017.66}{Flexible mechanical
  metamaterials}, Nature Reviews Materials 2~(11) (2017) 17066.

\bibitem{Kadic:2019aa}
M.~Kadic, G.~W. Milton, M.~van Hecke, M.~Wegener,
  \href{https://doi.org/10.1038/s42254-018-0018-y}{3d metamaterials}, Nature
  Reviews Physics 1~(3) (2019) 198--210.

\bibitem{10.1063/5.0152099}
J.~Zhang, B.~Hu, S.~Wang, \href{https://doi.org/10.1063/5.0152099}{{Review and
  perspective on acoustic metamaterials: From fundamentals to applications}},
  Applied Physics Letters 123~(1) (2023) 010502.

\bibitem{Geers2010}
M.~Geers, V.~Kouznetsova, W.~Brekelmans,
  \href{https://www.sciencedirect.com/science/article/pii/S0377042709005536}{Multi-scale
  computational homogenization: Trends and challenges}, Journal of
  Computational and Applied Mathematics 234~(7) (2010) 2175--2182, fourth
  International Conference on Advanced COmputational Methods in ENgineering
  (ACOMEN 2008).

\bibitem{hassani1998review}
B.~Hassani, E.~Hinton, A review of homogenization and topology optimization
  i---homogenization theory for media with periodic structure, Computers \&
  Structures 69~(6) (1998) 707--717.

\bibitem{Schmidt2022}
F.~Schmidt, M.~Kr{\"u}ger, M.-A. Keip, C.~Hesch, Computational homogenization
  of higher-order continua, International Journal for Numerical Methods in
  Engineering 123~(11) (2022) 2499--2529.

\bibitem{gautschi2002piezoelectric}
G.~Gautschi, Piezoelectric sensors, in: Piezoelectric Sensorics, Springer,
  2002, pp. 73--91.

\bibitem{sinha2009piezoelectric}
N.~Sinha, G.~E. Wabiszewski, R.~Mahameed, V.~V. Felmetsger, S.~M. Tanner, R.~W.
  Carpick, G.~Piazza, Piezoelectric aluminum nitride nanoelectromechanical
  actuators, Applied Physics Letters 95~(5) (2009) 053106.

\bibitem{safaei2019review}
M.~Safaei, H.~A. Sodano, S.~R. Anton, A review of energy harvesting using
  piezoelectric materials: state-of-the-art a decade later (2008--2018), Smart
  Materials and Structures 28~(11) (2019) 113001.

\bibitem{dagdeviren2016recent}
C.~Dagdeviren, P.~Joe, O.~L. Tuzman, K.-I. Park, K.~J. Lee, Y.~Shi, Y.~Huang,
  J.~A. Rogers, Recent progress in flexible and stretchable piezoelectric
  devices for mechanical energy harvesting, sensing and actuation, Extreme
  mechanics letters 9 (2016) 269--281.

\bibitem{guerin2021restriction}
S.~Guerin, D.~Thompson, Restriction boosts piezoelectricity, Nature Materials
  20~(5) (2021) 574--575.

\bibitem{jaffe1955properties}
B.~Jaffe, R.~Roth, S.~Marzullo, Properties of piezoelectric ceramics in the
  solid-solution series lead titanate-lead zirconate-lead oxide: tin oxide and
  lead titanate-lead hafnate, Journal of research of the National Bureau of
  Standards 55~(5) (1955) 239--254.

\bibitem{haertling1999ferroelectric}
G.~H. Haertling, Ferroelectric ceramics: history and technology, Journal of the
  American Ceramic Society 82~(4) (1999) 797--818.

\bibitem{jaffe1958piezoelectric}
H.~Jaffe, Piezoelectric ceramics, Journal of the American Ceramic Society
  41~(11) (1958) 494--498.

\bibitem{saito2004lead}
Y.~Saito, H.~Takao, T.~Tani, T.~Nonoyama, K.~Takatori, T.~Homma, T.~Nagaya,
  M.~Nakamura, Lead-free piezoceramics, Nature 432~(7013) (2004) 84--87.

\bibitem{EditorialJAP2022}
I.~Arias, G.~Catalan, P.~Sharma, The emancipation of flexoelectricity, Journal
  of Applied Physics 131~(020401) (2022).

\bibitem{zubko2013flexoelectric}
P.~Zubko, G.~Catalan, A.~K. Tagantsev, Flexoelectric effect in solids, Annual
  Review of Materials Research 43 (2013) 387--421.

\bibitem{Cross2006}
L.~E. Cross, \href{https://doi.org/10.1007/s10853-005-5916-6}{Flexoelectric
  effects: Charge separation in insulating solids subjected to elastic strain
  gradients}, Journal of Materials Science 41~(1) (2006) 53--63.

\bibitem{sharma2007possibility}
N.~Sharma, R.~Maranganti, P.~Sharma, On the possibility of piezoelectric
  nanocomposites without using piezoelectric materials, Journal of the
  Mechanics and Physics of Solids 55~(11) (2007) 2328--2350.

\bibitem{mocci2021geometrically}
A.~Mocci, J.~Barcel{\'o}-Mercader, D.~Codony, I.~Arias, Geometrically polarized
  architected dielectrics with apparent piezoelectricity, Journal of the
  Mechanics and Physics of Solids 157 (2021) 104643.

\bibitem{Mawassy2023}
N.~Mawassy, H.~Reda, J.-F. Ganghoffer, H.~Lakiss, Control of the piezoelectric
  and flexoelectric homogenized properties of architected materials by tuning
  their inner topology, Mechanics Research communication 127~(104034) (2023).

\bibitem{Greco2023}
F.~Greco, D.~Codony, H.~Mohammadi, S.~Fern{\'a}ndez-M{\'e}ndez, I.~Arias,
  \href{https://www.sciencedirect.com/science/article/pii/S0022509623002818}{Topology
  optimization of flexoelectric metamaterials with apparent piezoelectricity},
  Journal of the Mechanics and Physics of Solids (2023) 105477\href
  {https://doi.org/https://doi.org/10.1016/j.jmps.2023.105477}
  {\path{doi:https://doi.org/10.1016/j.jmps.2023.105477}}.

\bibitem{Abdollahi2014}
A.~Abdollahi, C.~Peco, D.~Mill\'an, M.~Arroyo, I.~Arias,
  \href{http://dx.doi.org/10.1063/1.4893974}{Computational evaluation of the
  flexoelectric effect in dielectric solids}, Journal of Applied Physics
  116~(9) (2014) 093502.

\bibitem{Abdollahi2015a}
A.~Abdollahi, D.~Mill\'an, C.~Peco, M.~Arroyo, I.~Arias,
  \href{http://dx.doi.org/10.1103/PhysRevB.91.104103}{Revisiting pyramid
  compression to quantify flexoelectricity: A three-dimensional simulation
  study}, Phys. Rev. B 91 (2015) 104103.

\bibitem{zhuang2020meshfree}
X.~Zhuang, S.~Nanthakumar, T.~Rabczuk, A meshfree formulation for large
  deformation analysis of flexoelectric structures accounting for the surface
  effects, Engineering Analysis with Boundary Elements 120 (2020) 153--165.

\bibitem{ghasemi2017level}
H.~Ghasemi, H.~S. Park, T.~Rabczuk, A level-set based iga formulation for
  topology optimization of flexoelectric materials, Computer Methods in Applied
  Mechanics and Engineering 313 (2017) 239--258.

\bibitem{codony2020modeling}
D.~Codony, P.~Gupta, O.~Marco, I.~Arias, Modeling flexoelectricity in soft
  dielectrics at finite deformation, Journal of the Mechanics and Physics of
  Solids 146 (2020) 104182.

\bibitem{codony2019immersed}
D.~Codony, O.~Marco, S.~Fern{\'a}ndez-M{\'e}ndez, I.~Arias, An immersed
  boundary hierarchical b-spline method for flexoelectricity, Computer Methods
  in Applied Mechanics and Engineering (2019).

\bibitem{codony2021mathematical}
D.~Codony, A.~Mocci, J.~Barcel{\'o}-Mercader, I.~Arias, Mathematical and
  computational modeling of flexoelectricity, Journal of Applied Physics
  130~(23) (2021) 231102.

\bibitem{mao2016mixed}
S.~Mao, P.~K. Purohit, N.~Aravas, Mixed finite-element formulations in
  piezoelectricity and flexoelectricity, Proceedings of the Royal Society A:
  Mathematical, Physical and Engineering Sciences 472~(2190) (2016) 20150879.

\bibitem{Deng2017}
F.~Deng, Q.~Deng, W.~Yu, S.~Shen, Mixed finite elements for flexoelectric
  solids, Journal of Applied Mechanics 84~(8) (2017) 081004.

\bibitem{Tian2021}
X.~Tian, J.~Sladek, V.~Sladek, Q.~Deng, Q.~Li, A collocation mixed finite
  element method for the analysis of flexoelectric solids, International
  Journal of Solids and Structures 217-218 (2021) 27--39.

\bibitem{ventura2020c0}
J.~Ventura, D.~Codony, S.~Fern{\'a}ndez-M{\'e}ndez,
  \href{https://doi.org/10.1007/s10915-021-01613-w}{A c0 interior penalty
  finite element method for flexoelectricity}, Journal of Scientific Computing
  88~(3) (2021) 88.

\bibitem{barcelo2022weak}
J.~Barcel{\'o}-Mercader, D.~Codony, S.~Fern{\'a}ndez-M{\'e}ndez, I.~Arias, Weak
  enforcement of interface continuity and generalized periodicity in high-order
  electromechanical problems, International Journal for Numerical Methods in
  Engineering 123~(4) (2022) 901--923.

\bibitem{Balcells}
O.~Balcells-Quintana, D.~Codony, S.~Fern{\'a}ndez-M{\'e}ndez, C0-ipm with
  generalised periodicity and application to flexoelectricity-based 2d
  metamaterials, Journal of Scientific Computing 313 (2022) 239--258.

\bibitem{hill1963elastic}
R.~Hill, Elastic properties of reinforced solids: some theoretical principles,
  Journal of the Mechanics and Physics of Solids 11~(5) (1963) 357--372.

\bibitem{hill1967essential}
R.~Hill, The essential structure of constitutive laws for metal composites and
  polycrystals, Journal of the Mechanics and Physics of Solids 15~(2) (1967)
  79--95.

\bibitem{deBoor2001}
C.~de~Boor, \href{http://www.springer.com/gb/book/9780387953663}{A Practical
  Guide to Splines}, Applied Mathematical Sciences, Springer New York, 2001.
\newline\urlprefix\url{http://www.springer.com/gb/book/9780387953663}

\bibitem{burman2010ghost}
E.~Burman, Ghost penalty, Comptes Rendus Mathematique 348~(21-22) (2010)
  1217--1220.

\bibitem{Duster2008}
A.~D\"{u}ster, J.~Parvizian, Z.~Yang, E.~Rank,
  \href{http://dx.doi.org/10.1016/j.cma.2008.02.036}{The finite cell method for
  three-dimensional problems of solid mechanics}, Computer Methods in Applied
  Mechanics and Engineering 197~(45) (2008) 3768 -- 3782.

\bibitem{schillinger2015finite}
D.~Schillinger, M.~Ruess, The finite cell method: A review in the context of
  higher-order structural analysis of cad and image-based geometric models,
  Archives of Computational Methods in Engineering 22~(3) (2015) 391--455.

\bibitem{Sharma2010}
N.~Sharma, C.~Landis, P.~Sharma,
  \href{http://dx.doi.org/10.1063/1.3443404}{Piezoelectric thin-film
  superlattices without using piezoelectric materials}, Journal of Applied
  Physics 108~(2) (2010) 1--25.

\bibitem{Mindlin1968a}
R.~D. Mindlin, N.~N. Eshel,
  \href{http://dx.doi.org/10.1016/0020-7683(68)90036-X}{On first
  strain-gradient theories in linear elasticity}, International Journal of
  Solids and Structures 4~(1) (1968) 109--124.

\bibitem{Altan1997}
B.~S. Altan, E.~C. Aifantis, On some aspects in the special theory of gradient
  elasticity, Journal of the Mechanical Behavior of Materials 8~(3) (1997)
  231--282.

\end{thebibliography}
\end{document}